\documentclass[proc]{edpsmath}

\usepackage[utf8]{inputenc}
\usepackage[T1]{fontenc}
\usepackage{lmodern}

\usepackage{amsfonts,amssymb,amsmath,amsthm}
\usepackage{booktabs}
\usepackage{color}
\usepackage{cite}
\usepackage{caption}
\usepackage{fullpage}
\usepackage{graphicx,hyperref}
\usepackage{mathabx,mathrsfs,mathtools}
\usepackage{verbatim}
\usepackage{cancel}
\usepackage{tikz,pgfplots}
\usepackage{tabu}
\usepackage{subcaption}
\usepackage[active]{srcltx}
\usepackage[official]{eurosym}
\usepackage[noend]{algpseudocode}
\usepackage{algorithm}
\usepackage{accents}
\usepackage[yyyymmdd,hhmmss]{datetime}

\allowdisplaybreaks
\hypersetup{
  linkcolor  = green!60!black,
  citecolor  = blue!60!black,
  urlcolor   = red!60!black,
  colorlinks = true,
}

\numberwithin{equation}{section}
\definecolor{jpred}{rgb}{0.6,0.2,0.2}

\def\argmax{\mathop{\rm argmax}}
\def\argmin{\mathop{\rm argmin}}

\makeatletter

\numberwithin{dummy}{section}

\newtheorem{assumption}{Assumption}[section]
\newtheorem{remark}{Remark}[section]

\newtheorem*{conjecture*}{Conjecture}
\makeatother

\def\argmin_#1{\underset{#1}{\mathrm{argmin\, }}}
\def\argmax_#1{\underset{#1}{\mathrm{argmax\, }}}

\def \Sup{\displaystyle\sup}

\def \trans{^{\scriptscriptstyle{\intercal}}}

\def \Sup{\displaystyle\sup}

\def \Max{\displaystyle\max}

\def \b1{\bf{1}}

\def\beqs{\begin{eqnarray*}}
\def\enqs{\end{eqnarray*}}
\def\beq{\begin{eqnarray}}
\def\enq{\end{eqnarray}}

\def\bk{{\bf k}}

\def \N{\mathbb{N}}
\def \R{\mathbb{R}}
\def \L{\mathbb{L}}
\def \M{\mathbb{M}}

\def \E{\mathbb{E}}
\def \F{\mathbb{F}}

\def \P{\mathbb{P}}

\def \Ac{{\cal A}}
\def \Bc{{\cal B}}
\def \Cc{{\cal C}}

\def \Fc{{\cal F}}

\def \Lc{{\cal L}}
\def \Pc{{\cal P}}

\def \Nc{{\cal N}}

\def \Wc{{\cal W}}

\def \eps{\varepsilon}

\def\reff#1{{\rm(\ref{#1})}}

\def\beqs{\begin{eqnarray*}}
  \def\enqs{\end{eqnarray*}}
\def\beq{\begin{eqnarray}}
  \def\enq{\end{eqnarray}}

\begin{document}

\title{A class of finite-dimensional numerically solvable McKean-Vlasov control problems}
\runningtitle{}
\thanks{\noindent}
\author{Alessandro Balata}\address{University of Leeds, Leeds, United Kingdom;\\
\email{A.Balata@leeds.ac.uk}}
\author{C\^ome Hur\'e}\address{Univ. Paris Diderot - LPSM, Paris, France;\\
\email{hure@lpsm.paris}}
\author{Mathieu Lauri\`ere}\address{Operations Research and Financial Engineering, Princeton University, Princeton, USA;\\
\email{lauriere@princeton.edu}}
\author{Huy\^en Pham}\address{Univ. Paris Diderot - LPSM, Paris, France;\\
\email{pham@math.univ-paris-diderot.fr}}
\author{Isaque Pimentel}\address{CMAP, Ecole Polytechnique, Palaiseau, France;\\
\email{isaque.santa-brigida-pimentel@polytechnique.edu}}
\begin{abstract}
We address a class of McKean-Vlasov (MKV) control problems with common noise, called polynomial conditional MKV, and extending the known class of linear quadratic
stochastic MKV control problems.
We show how  this polynomial  class can be reduced by suitable Markov embedding to finite-dimensional stochastic control problems, and provide a
discussion and comparison of three probabilistic numerical methods for solving the reduced control problem: quantization, regression by control randomization, and regress-later methods.
Our numerical results are illustrated on various examples from portfolio selection and liquidation under drift uncertainty, and a model of interbank systemic risk with partial observation.
\end{abstract}

\maketitle

\vspace{2mm}



\newpage

\section{Introduction}\label{intro}

The optimal control of  McKean-Vlasov (also called mean-field) dynamics is a rather new topic in the area of stochastic control and applied probability, which has been knowing a surge of interest with the emergence of the mean-field game theory.
It is motivated on the one hand by the asymptotic formulation of cooperative equilibrium for a large population of particles (players) in mean-field interaction, and on the other hand from control problems with cost functional  involving nonlinear functional of the law of the state process (e.g.,~the mean-variance portfolio selection problem or risk measure  in finance).

In this paper, we are interested in McKean-Vlasov (MKV) control problems under partial observation and common noise, whose formulation is described as follows.
On a probability space $(\Omega,\Fc,\P)$ equipped with two independent Brownian motions $B$ and $W^0$,  let us consider the controlled stochastic MKV
dynamics in $\R^n$:
\begin{equation}
\label{dynX}
dX_s = b(X_s,\P_{_{X_s}}^{W^0},\alpha_s) ds +  \sigma(X_s,\P_{_{X_s}}^{W^0},\alpha_s) dB_s +  \sigma_0(X_s,\P_{_{X_s}}^{W^0},\alpha_s) dW_s^0, \quad X_0=x_0 \in \R^n,
\end{equation}
where $\P_{_{X_s}}^{W^0}$ denotes the conditional distribution of $X_s$ given $W^0$ (or equivalently given $\Fc_s^0$ where $\F^0$ $=$ $(\Fc_t^0)_t$ is the natural filtration generated by $W^0$), and the control $\alpha$ is $\F^0$-progressive valued in some Polish space $A$.
This measurability condition on the control means that the controller has a partial observation of the state, in the sense that she can only observe the common noise.
We make the standard Lipschitz condition on the coefficients $b(x,\mu,a), \; \sigma(x,\mu,a), \; \sigma_0(x,\mu,a) $ with respect to $(x,\mu)$ in
$\R^n\times \Pc_{_2}(\R^n)$, uniformly in $a$ $\in$ $A$, where  $\Pc_{_2}(\R^n)$
is  the set of all probability measures on $(\R^n,\Bc(\R^n))$ with a finite second-order moment, endowed with the $2$-Wasserstein metric $\Wc_{_2}$.
This ensures the well-posedness of the controlled MKV
stochastic differential equation (SDE) (\ref{dynX}).
The cost functional over a finite  horizon $T$ associated to the stochastic MKV
equation \reff{dynX}  (sometimes called conditional MKV
equation) for a control process $\alpha$,  is
\beqs
J(\alpha) &=& \E \Big[ \int_0^T f(X_t,\P_{_{X_t}}^{W^0},\alpha_t) dt + g(X_T,\P_{_{X_T}}^{W^0}) \Big],
\enqs
and the objective is to maximize over an admissible set $\Ac$ of control processes the cost functional:
\beq \label{controlMKV}
V_0 &=& \Sup_{\alpha\in\Ac} J(\alpha).
\enq
The set $\Ac$ of admissible controls usually requires some integrability conditions depending on the growth conditions on $f$, $g$, in order to ensure that $J(\alpha)$ is well-defined for
$\alpha$ $\in$ $\Ac$ (more details will be given in the examples, see Section~\ref{sec:numres}).
Notice that classical partial observation control problem (without
MKV dependence on the coefficients) arises as a particular case of
\reff{dynX}-\reff{controlMKV}.
We refer to the introduction in \cite{phawei17} for the details.

Let us recall from \cite{phawei17} the dynamic programming equation associated to the conditional
MKV control problem \reff{controlMKV}.
We start by defining a suitable dynamic version of this problem.
Let us consider $\Fc_0$ a sub $\sigma$-algebra of $\Fc$ independent of $B,W^0$.
It is assumed w.l.o.g. that $\Fc_0$ is rich enough in the sense that $\Pc_{_2}(\R^n)$ $=$
$\{\Lc(\xi): \xi \in L^2(\Fc_0;\R^n)\}$, where $\Lc(\xi)$ denotes the law of $\xi$.
Given a control $\alpha$ $\in$ $\Ac$, we consider the dynamic version of \reff{dynX} starting from $\xi$ $\in$ $L^2(\Fc_0;\R^n)$ at time $t$ $\in$ $[0,T]$, and written as:
\beqs
X_s^{t,\xi,\alpha} &=& \xi + \int_t^s b(X_u^{t,\xi,\alpha},\P_{_{X_u^{t,\xi,\alpha}}}^{W^0},\alpha_u) du
+ \int_t^s \sigma(X_u^{t,\xi,\alpha},\P_{_{X_u^{t,\xi,\alpha}}}^{W^0},\alpha_u) dB_u  \\
& & \;\;\; + \;  \int_t^s \sigma_0(X_u^{t,\xi,\alpha},\P_{_{X_u^{t,\xi,\alpha}}}^{W^0},\alpha_u) dW_u^0  \;,\;   t \leq s \leq  T.
\enqs
Let us then define the dynamic cost functional:
\beqs
J(t,\xi,\alpha) &=& \E \Big[ \int_t^T f(X_s^{t,\xi,\alpha},\P_{_{X_s^{t,\xi,\alpha}}}^{W^0},\alpha_s) ds + g(X_T^{t,\xi,\alpha},\P_{_{X_T^{t,\xi,\alpha}}}^{W^0}) \Big],
\enqs
for $(t,\xi) \in [0,T]\times L^2(\Fc_0;\R^n), \alpha\in\Ac$, and notice by the law of conditional expectations, and as $\alpha$ is $\F^0$-progressive that
\beqs
J(t,\xi,\alpha) &=& \E \Big[ \int_t^T \hat f(\P_{_{X_s^{t,\xi,\alpha}}}^{W^0},\alpha_s) ds + \hat g(\P_{_{X_s^{t,\xi,\alpha}}}^{W^0}) \Big],
\enqs
where $\hat f$ $:$ $\Pc_{_2}(\R^n)\times A$ $\rightarrow$ $\R$, $\hat g$ $:$ $\Pc_{_2}(\R^n)$ $\rightarrow$ $\R$ are defined by
\beq
\hat f(\mu,a) \; = \;  \mu (f(\cdot,\mu,a)) & = &  \int_{\R^n} f(x,\mu,a) \mu(dx),  \label{defhatf} \\
\hat g(\mu) \; = \; \mu(g(\cdot,\mu)) &=&    \int_{\R^n} g(x,\mu) \mu(dx).   \label{defhatg}
\enq
Moreover, notice that the conditional  law of $X_s^{t,\xi,\alpha}$ given $W^0$ depends on $\xi$ only through its law $\Lc(\xi)$, and we can then define for $\alpha$ $\in$
$\Ac$:
\beqs
\rho_s^{t,\mu,\alpha} &=& \P_{_{X_s^{t,\xi,\alpha}}}^{W^0}, \;\;\; \mbox{ for } \;  t \leq s, \;  \mu = \Lc(\xi)  \in \Pc_{_2}(\R^n).
\enqs
Therefore, the dynamic cost functional $J(t,\xi,\alpha)$ depends on $\xi$ $\in$ $L^2(\Fc_0;\R^n)$ only through its law $\Lc(\xi)$, and by an abuse of notation, we write
$J(t,\mu,\alpha)$ $=$ $J(t,\xi,\alpha)$ when $\mu$ $=$ $\Lc(\xi)$.
We then consider the value function for the conditional
MKV control problem \reff{controlMKV}, defined on
$[0,T]\times\Pc_{_2}(\R^n)$ by
\beq \label{defv}
v(t,\mu) &  = & \Sup_{\alpha\in\Ac} J(t,\mu,\alpha)
\; = \; \Sup_{\alpha\in\Ac} \E \Big[ \int_t^T \hat f(\rho_s^{t,\mu,\alpha},\alpha_s) ds + \hat g(\rho_T^{t,\mu,\alpha}) \Big],
\enq
and notice that at time $t$ $=$ $0$, when $\xi$ $=$ $x_0$ is a constant, then $V_0$ $=$ $v(0,\delta_{x_0})$.

It is shown in \cite{phawei17} that dynamic programming principle (DPP) for the conditional
MKV control problem \reff{defv} holds:  for $(t,\mu) \in [0,T]\times\Pc_{_2}(\R^n)$,
\beqs
v(t,\mu) &=& \Sup_{\alpha\in\Ac} \E \Big[  \int_t^\theta \hat f(\rho_s^{t,\mu,\alpha},\alpha_s) ds + v(\theta,\rho_\theta^{t,\mu,\alpha}) \Big],
\enqs
for any $\F^0$-stopping time $\theta$ valued in $[t,T]$.
Next, by relying on the notion of differentiability with respect to probability measures introduced by P. L. Lions \cite{lio12} (see also the lecture notes \cite{car12}) and the chain rule (It\^o's formula) along flow of probability measures (see \cite{buetal14}, \cite{chacridel15}), we derive the
HJB equation for $v$:
\begin{equation} \label{HJBdynpro}
\left\{
\begin{array}{rcl}
\partial_t v  +  \Sup_{a\in A} \Big[ \hat f(\mu,a) + \mu\big(\L^a v(t,\mu) \big) + \mu\otimes\mu\big(\M^a v(t,\mu) \big) \Big] & =& 0,  \;\;\; (t,\mu)  \in [0,T)\times\Pc_{_2}(\R^n), \\
v(T,\mu) &=& \hat g(\mu), \;\;\; \mu \in \Pc_{_2}(\R^n),
\end{array}
\right.
\end{equation}
where for $\phi$ $\in$ $\Cc_b^2(\Pc_{_2}(\R^n))$, $a$ $\in$ $A$, and $\mu$ $\in$ $\Pc_{_2}(\R^n)$, $\L^a\phi(\mu)$ is the function $\R^n$ $\rightarrow$ $\R$ defined by
\beq \label{defL}
\L^a\phi(\mu)(x) &=& \partial_\mu \phi(\mu)(x).b(x,\mu,a) + \frac{1}{2}{\rm tr}\big(\partial_x\partial_\mu\phi(\mu)(x)( \sigma\sigma\trans + \sigma_0\sigma_0\trans)(x,\mu,a) \big),
\enq
and $\M^a\phi(\mu)$  is the function $\R^n\times\R^n$ $\rightarrow$ $\R$ defined by
\beq \label{defM}
\M^a\phi(\mu)(x,x') &=& \frac{1}{2} {\rm tr}\big(  \partial_\mu^2\phi(\mu)(x,x')\sigma_0(x,\mu,a)\sigma_0\trans(x',\mu,a) \big).
\enq

The HJB equation \reff{HJBdynpro} is a fully nonlinear partial differential equation (PDE) in the infinite-dimensional Wasserstein space.
In general, this PDE does not have an explicit solution except in the notable important class of linear-quadratic
MKV control problem.
Numerical resolution for MKV control problem or equivalently for the associated HJB equation is a challenging problem due to the nonlinearity of the optimization problem and the infinite-dimensional feature of the Wasserstein space.
In this work, our purpose is to investigate a class of MKV control problems which can be reduced to finite-dimensional problems in view of numerical resolution.

\section{Polynomial McKean-Vlasov control problem}

\subsection{Main assumptions}

We make two kinds of assumptions on the coefficients of the model: one on the dependence on $x$ and the other on the dependence on $\mu$.

\vspace{1mm}

\noindent {\bf Assumptions: dependence on $x$:}
we consider a class of models where the coefficients of the MKV equation are linear w.r.t. the state variable $X$, i.e., they are in the form
\begin{equation} \label{dynblin}
\left\{
\begin{array}{ccc}
b(x,\mu,a) & = & b_0(\mu,a) + b_1(\mu,a)x, \\
\sigma(x,\mu,a) &=& \vartheta_0(\mu,a) + \vartheta_1(\mu,a) x, \\
\sigma_0(x,\mu,a) &=& \gamma_0(\mu,a) + \gamma_1(\mu,a) x,
\end{array}
\right.
\end{equation}
while the running and terminal cost functions are polynomial in the state variable in the sense that they are in the form
\beqs
f(x,\mu,a) &=& f_0(\mu,a) + f_1(\mu,a) x +  \sum_{k=2}^p f_k(\mu,a) |x|^k, \\
g(x,\mu) &=& g_0(\mu) + g_1(\mu) x +  \sum_{k=2}^p g_k(\mu) |x|^k,
\enqs
for some integer $p$ $\geq$ $2$ (here, $f_1$ and $g_1$ are functions taking value in $\R^n$, and $|\cdot|$ denotes the Euclidean norm).

\vspace{1mm}

\noindent {\bf Assumptions: dependence on $\mu$:} we  assume that all the coefficients depend on $\mu$ through its first $p$ moments, i.e., they are in the form
\begin{equation} \label{formbf}
\left\{
\begin{array}{ccc}
b_0(\mu,a) \; = \; \bar b_0(\bar\mu,\bar\mu_2,\ldots,\bar\mu_p,a), & & b_1(\mu,a) \; = \; \bar b_1(\bar\mu,\bar\mu_2,\ldots,\bar\mu_p,a) \\
\vartheta_0(\mu,a) \; = \; \bar\vartheta_0(\bar\mu,\bar\mu_2,\ldots,\bar\mu_p,a), & & \vartheta_1(\mu,a) \; = \; \bar\vartheta_1(\bar\mu,\bar\mu_2,\ldots,\bar\mu_p,a) \\
\gamma_0(\mu,a) \; = \; \bar\gamma_0(\bar\mu,\bar\mu_2,\ldots,\bar\mu_p,a), & & \gamma_1(\mu,a) \; = \; \bar\gamma_1(\bar\mu,\bar\mu_2,\ldots,\bar\mu_p,a) \\
f_k(\mu,a) \; = \; \bar f_k(\bar\mu,\bar\mu_2,\ldots,\bar\mu_p,a), & & g_k(\mu) \; = \; \bar g_k(\bar\mu,\bar\mu_2,\ldots,\bar\mu_p), \;\;\;\;\;  k=0,\ldots,p,
\end{array}
\right.
\end{equation}
where, given $\mu$ $\in$ $\Pc_{p}(\R^n)$,  we denote by
\beqs
\bar\mu \; = \;  \int_{\R^n} x \mu(dx), & &  \bar\mu_k \; = \; \int_{\R^n} |x|^k  \mu(dx), \;\;\;\; k=2,\ldots,p.
\enqs

 \noindent We assume that the coefficients $\bar b_0, \bar b_1, \bar \vartheta_0, \bar \vartheta_1, \bar \gamma_0, \bar \gamma_1$ are Lipschitz w.r.t. the $p$ first arguments uniformly w.r.t. the control argument $a \in A$.
This condition will ensure existence and uniqueness of a solution to the finite-dimensional
MKV SDE defined later in \eqref{dynXreduc}.

Notice that in this case, the functions $\hat f$ and $\hat g$ defined in \reff{defhatf}-\reff{defhatg} are given by
 \beqs
 \hat f(\mu,a) &=& \bar f_0(\bar\mu,\bar\mu_2,\ldots,\bar\mu_p,a) + \bar f_1(\bar\mu,\bar\mu_2,\ldots,\bar\mu_p,a) \bar\mu
 +\sum_{k=2}^p \bar f_k(\bar\mu,\bar\mu_2,\ldots,\bar\mu_p,a) \bar\mu_k \\
 &=:& \bar f(\bar\mu,\bar\mu_2,\ldots,\bar\mu_p,a), \\
 \hat g(\mu) &=&  \bar g_0(\bar\mu,\bar\mu_2,\ldots,\bar\mu_p) +  \bar g_1(\bar\mu,\bar\mu_2,\ldots,\bar\mu_p) \bar\mu
 + \sum_{k=2}^p \bar g_k(\bar\mu,\bar\mu_2,\ldots,\bar\mu_p) \bar\mu_k \\
 &=:& \bar g(\bar\mu,\bar\mu_2,\ldots,\bar\mu_p).
 \enqs

\begin{remark}
 A more general class of running and terminal cost functions would be to consider multi-polynomial of degree $p$ functions $f$ and $g$ in the form
\beqs
f(x,\mu,a) \; = \;  \sum_{|{\bf k}|=0}^p f_{\bk}\left((\mu^{\bk'})_{|\bk'| \leq p},a\right)  x^{\bk}, & &  g(x,\mu) \; = \; \sum_{|{\bf k}|=0}^p g_{\bk}\left((\mu^{\bk'})_{|\bk'| \leq p}\right)  x^{\bk},
\enqs
where we use multi-index notations $\bk$ $=$ $(k_1,\ldots,k_n)$ $\in$ $\N^n$, $|\bk|$ $=$ $k_1+\ldots+k_n$, $x^{\bk}$ $=$ $x_1^{k_1}\ldots x_n^{k_n}$ for $x$ $=$
$(x_1,\ldots,x_n)$ $\in$ $\R^n$
and
\beqs
\mu^{\bk} &=& \int_{\R^n}  x^{\bk} \mu(dx).
\enqs
\end{remark}

\subsection{Markovian embedding}

Given the controlled process $X$ $=$ $X^\alpha$ solution to the stochastic
MKV dynamics \reff{dynX}, denote by
\beqs
\bar X_t \; = \;  \E[ X_t | W^0], & & Y_t^k \; = \; \E[ |X_t|^k | W^0], \;\;\; k =2, \ldots,p.
\enqs
To alleviate the notations, let us assume that $n=1$ (otherwise multi-indices should be used).
From the linear/polynomial assumptions \reff{dynblin}-\reff{formbf}, by It\^o's formula and taking conditional expectations, we can derive the dynamics of $(\bar X,Y^2,\ldots,Y^p)$ as
\begin{equation} \label{dynXreduc}
\left\{
\begin{array}{ccl}
d \bar X_t &=&  \bar B(\bar X_t,Y_t^2,\ldots,Y_t^p,\alpha_t) dt +  \bar\Sigma(\bar X_t,Y_t^2,\ldots,Y_t^p,\alpha_t) dW_t^0,  \quad \bar{X}_0=x_0,\\
dY_t^k &=&  B_k(\bar X_t,Y_t^2,\ldots,Y_t^p,\alpha_t) dt  +  \Sigma_k(\bar X_t,Y_t^2,\ldots,Y_t^p,\alpha_t) dW_t^0, \quad Y^k_0= \vert x_0 \vert^k, \; \;\; k=2,\ldots,p,
\end{array}
\right.
\end{equation}
where
\beqs
\bar B(\bar x,y^2,\ldots,y^p,a) &=&  \bar b_0(\bar x,y^2,\ldots,y^p,a) +  \bar b_1(\bar x,y^2,\ldots,y^p,a) \bar x,
\\
\bar\Sigma(\bar x,y^2,\ldots,y^p,a) &=&  \bar\gamma_0(\bar x,y^2,\ldots,y^p,a) +  \bar\gamma_1(\bar x,y^2,\ldots,y^p,a) \bar x,
\\
B_k(\bar x,y^2,\ldots,y^p,a) &=&  k \bar b_0(y^2,\dots,y^p,a) y^{k-1} + k \bar b_1(y^2,\dots,y^p,a) y^{k}
\\
  &&+\frac{k(k-1)}{2}(\bar \vartheta_0(y^2,\dots,y^p,a))^2 y^{k-2} + \frac{k(k-1)}{2}(\bar \vartheta_1(y^2,\dots,y^p,a))^2 y^k
  \\
  &&+ k(k-1) \bar \vartheta_0(y^2,\dots,y^p,a) \bar \vartheta_1(y^2,\dots,y^p,a) y^{k-1}
  \\
  &&+\frac{k(k-1)}{2}(\bar \gamma_0(y^2,\dots,y^p,a))^2 y^{k-2} +\frac{k(k-1)}{2} (\bar \gamma_1(y^2,\dots,y^p,a))^2 y^k
  \\
  &&+ k(k-1) \bar \gamma_0(y^2,\dots,y^p,a) \bar \gamma_1(y^2,\dots,y^p,a) y^{k-1},
  \\
\Sigma_k(\bar x,y^2,\ldots,y^p,a) &=& k \left( \bar \gamma_0(y^2,\dots,y^p,a) y^{k-1} + \bar \gamma_1(y^2,\dots,y^p,a) y^k \right),
\enqs
while the cost functional is written as
\beq \label{Jreduc}
J(\alpha) &=& \E \Big[ \int_0^T  \bar f(\bar X_t, Y_t^2,\ldots,Y_t^p,\alpha_t) dt +   \bar g(\bar X_T, Y_T^2,\ldots,Y_T^p) \Big].
\enq

The
MKV control problem is then reduced in this polynomial framework into a finite-dimensional control problem with $\F^0$-adapted controlled variables  $(\bar X,Y^2,\ldots,Y^p)$.
In the next section, we describe three probabilistic numerical methods for solving finite-dimensional stochastic control problems and will apply in section~\ref{sec:numres} each of these methods to three examples arising from polynomial MKV control problems under partial observation and common noise.

\newcommand*{\backin}{\rotatebox[origin=c]{-180}{$\in$}}

\newlength{\dhatheight}
\newcommand{\doublehat}[1]{
  \settoheight{\dhatheight}{\ensuremath{\hat{#1}}}
  \addtolength{\dhatheight}{-0.35ex}
  \hat{\vphantom{\rule{1pt}{\dhatheight}}
    \smash{\hat{#1}}}}

\section{Numerical methods}

In this section, we introduce our numerical methods for the resolution of the reduced problem \reff{dynXreduc}-\reff{Jreduc}.

Let us introduce the process $Z^\alpha$, valued in $\R^d$, controlled by an adapted process $\alpha$ taking values in $A$, solution to
\begin{equation}
  d Z_t^\alpha = b(Z_t^\alpha, \alpha_t) dt + \sigma_0(Z_t^\alpha, \alpha_t) d W_t^0, \quad Z_0^\alpha= z_0 \in \R^d,
\end{equation}
and the performance measure
\begin{equation}
  J(t,z,\alpha) = \E\left[ \int_t^T f(Z_t^\alpha, \alpha_t) dt + g(Z_T^\alpha) \Big\vert Z^\alpha_t=z\right],
\end{equation}
which assesses the average performance of the control.

Introduce now a time discretization  $t_n=n\Delta t$, $n$ $=$ $0,\ldots,N$, $\Delta t$ $=$ $T/N$, and denote by $\mathcal{A}_{\Delta t}$ the space of discrete processes $(\alpha_{t_n})_{n=0}^{N-1}$ such that for all $n$, $n=0,\ldots,N-1$, $\alpha_{t_n}$ is $\mathcal{F}^0_{t_n}$-measurable.

We can write the Euler approximation of the SDE governing the process $Z$ $=$ $Z^\alpha$,  with $\alpha \in \mathcal{A}_{\Delta t}$ (to alleviate notations, we sometimes omit the dependence on $\alpha$ when there is no ambiguity, and keep the same notation $Z$ for the discrete and continuous process)
\begin{equation}
\label{eq:euler}
Z_{t_{n+1}} =Z_{t_{n}} + b(Z_{t_{n}}, \alpha_{t_{n}}) \Delta t + \sigma_0(Z_{t_{n}}, \alpha_{t_{n}}) \Delta W_{t_{n}}^0,
\end{equation}
where $\Delta W_{t_{n}}^0 \sim \mathcal{N}(0,\Delta t)$ is an increment of $W^0$.

The discrete time approximation of $J(t_n,z,\alpha)$ is defined as:
\begin{equation}
\label{eq:performance_measure}
J_{\Delta t}(t_{n},z,\alpha) = \E\left[ \sum_{k=n}^{N-1} f(Z_{t_{k}}, \alpha_{t_{k}})\Delta t + g(Z_{t_{N}}) \Big\vert Z_{t_n}=z\right],
\end{equation}
where $\alpha \in \mathcal{A}_{\Delta t}$ .

\subsection{Value and Performance iteration}

For $n=0,\ldots,N$, consider $V_{\Delta t }(t_{n},z)=\sup\limits_{\alpha \in \mathcal{A}_{\Delta t}} J_{\Delta t}(t_{n},z,\alpha)$, the discrete time approximation of the value function at time $t_n$: ${V\big (t_n,z\big)=\sup\limits_{\alpha \in \mathcal{A}} J(t_{n},z,\alpha)}$.
The dynamic programming principle states that $\big( V_{\Delta t}(t_n,\cdot) \big)_{0\leq n\leq N}$ is solution to the Bellman equation:
\begin{equation}
\label{eq:DPE}
\begin{cases}
& V^{}_{\Delta t}(t_N,z)=g(z)\\
&V^{}_{\Delta t}(t_n,z)=\Sup_{a \in A}\Big\{ f(z,a)\Delta t+\E_{n,z}^a\big[V^{}_{\Delta t}(t_{n+1},Z_{t_{n+1}}) \big] \Big\}, \;\; n = N-1,\ldots,0,
\end{cases}
\end{equation}
where $\E_{n,z}^a[\cdot]$ denotes the expectation conditioned on the event $\{Z_{t_n}=z\}$ and when using the control $\alpha_{t_n}$ $=$ $a$ at time $t_n$.
Observe that for $n=0,\ldots,N-1$, the equation \eqref{eq:DPE} provides a backward procedure to recursively compute  the $V^{}_{\Delta t}(t_n,\cdot)$ if we know how to analytically compute the conditional expectations $\E_{n,z}^a[V^{}_{\Delta t}(t_{n+1},Z_{t_{n+1}})]$ for all $z \in \R^d$ and all control $a\in A$.
We refer to the procedure in \eqref{eq:DPE} as value iteration.

An alternative approach to compute $V_{\Delta t}(t_n,\cdot)$, for $n=0,\ldots,N-1$, is to notice that once again by the dynamic programming principle, it holds that $\big( V_{\Delta t}(t_n,\cdot) \big)_{0\leq n\leq N}$ is solution to the backward equation
\begin{equation}
\label{eq:DPEpi}
\begin{cases}
V^{}_{\Delta t}(t_N,z)=g(z)\\
V^{}_{\Delta t}(t_n,z) =  \Sup_{a \in A} \; \left\{   f(z, a) \Delta t + \E_{n,z}^a \left[ \sum_{k=n+1}^{N-1} f\big(Z_{t_k}, \alpha^*_{t_k}(Z_{t_k})\big)  \Delta t+g(Z_{t_N}) \right] \right\},  \;\; n = N-1,\ldots,0,
\end{cases}
\end{equation}
where for $k=n+1,\ldots,N-1$,  the control $\alpha^*_{t_k}$ is the optimal control at time $t_k$ defined as follows:
\begin{equation}
\label{eq:control}
\alpha^*_{t_k}(z)  = \argmax_{a\in A}\left\{ f(z,a)\Delta t+ \E_{k,z}^a \left[ \sum_{\ell =k+1}^{N-1} f\big(Z^*_{t_\ell},\alpha^*_{t_\ell }(Z^*_{t_\ell })\big) \Delta t+g(Z^*_{t_N}) \right]  \right\},
\end{equation}
and where $\big(Z^*_{t_k}\big)_{n \leq k \leq N}$ is the process $Z$ controlled by the following control $\alpha$ from time $t_n$ to $t_N$:
\begin{equation}
	\begin{cases}
	\alpha_{t_{n}}=a,\\
	\alpha_{t_{k}}= \alpha^*_{t_{k}} \text{ for } n+1 \leq k \leq N-1.
	\end{cases}
\end{equation}

For $n=0,\ldots,N-1$, the scheme \eqref{eq:DPEpi} provides once again a backward procedure to compute $V^{}_{\Delta t}(t_n,\cdot)$, assuming that we know how to analytically compute the conditional expectations $\E_{n,z}^a \left[ \sum_{k=n+1}^{N-1} f\big(Z_{t_k}, \alpha^*_{t_k}(Z_{t_k})\big) \Delta t+g(Z_{t_N}) \right]$ for all $z \in \R^d$ and all control $a\in A$.
We refer to the procedure in \eqref{eq:DPEpi} as performance iteration. \\

\noindent Except for trivial cases, closed-form formulas for the conditional expectations appearing in the value and performance iteration procedures are not available, and they have to be approximated, which is the main difficulty when implementing both approaches to compute the value functions.
In the next section, we discuss different ways to approximate conditional expectations and derive the corresponding estimations of the value functions $V^{}_{\Delta t}(t_n,\cdot)$ for $n=0,\ldots,N-1$.

\vspace{3mm}
\subsection{Approximation of conditional expectations}
\label{sec:approxexpect}
In this subsection, we present three numerical methods that we apply later to conditional MKV problems.
Two of these methods belong to the class of Regression Monte Carlo techniques, a family of algorithms whose effectiveness highly relies on the choice of the basis functions used to approximate conditional expectations; the third algorithm, Quantization, approximate the controlled process $Z_{t_{n}}^\alpha$ with a particular finite state Markov chain for which expectations can be approximated quickly.

\subsubsection{Regression Monte Carlo}

In the simpler uncontrolled case, the family of Regression Monte Carlo algorithms is based on the idea of approximating the conditional expectation $\E\big[V_{\Delta t}(t_{n+1},Z_{t_{n+1}}) \big \vert Z_{t_n}\big] $, for $n=0,\ldots,N-1$, by the orthogonal projection of $V_{\Delta t}(t_{n+1},Z_{t_{n+1}})$ onto the space generated by a finite family of $\big\{ \phi_k(Z_{t_n}) \big\}_{k\geq 1}$ where $(\phi_k)_{k\geq 1}$ is a family of \emph{basis functions}, i.e., a family of measurable real-valued functions defined on $\R^d$ such that  $\big(\phi_k(Z_{t_n})\big)_{k\geq 1}$ is total in $L^2(\sigma(Z_{t_n}))$\footnote{$L^2(\sigma(Z_{t_n}))$ is the space of the square-integrable $\sigma(Z_{t_n})$-measurable r.v.} and such that for all scalars $\beta_{k}$ and all $K\geq 1$, if $\sum_{k=1}^{K} \beta_k \phi_k(Z_{t_n}) =0$ a.s. then $\beta_k=0, \; \text{ for } k=1,\ldots, K$.

The expectation $\E\big[V_{\Delta t}(t_{n+1},Z_{t_{n+1}}) \big \vert Z_{t_n}\big] $ should then be approximated as follows:
\begin{equation}
  \E\big[V_{\Delta t}(t_{n+1},Z_{t_{n+1}}) \big \vert Z_{t_n}\big] \approx \; \sum_{k=1}^{K} \beta^n_k \phi_{k}(Z_{t_n}),
\end{equation}
where $K\geq 1$ is fixed and $ \beta^{n} = ( \beta^{n}_{1},\,\dots,\, \beta^{n}_{K})^T$ is defined as:
\begin{equation}
\label{eq:regression_coeff_proj_old}
\beta^{n} \; = \; \argmin_{\beta \in \mathbb{R}^K }\left\{\E \left[ \Big| V_{\Delta t}(t_{n+1},Z_{t_{n+1}}) - \sum_{k=1}^K\beta_k\phi_{k}(Z_{t_{n}})  \Big|^2 \right] \right\}.
\end{equation}
Notice that $\beta^{n}$ is defined in \eqref{eq:regression_coeff_proj_old} as the minimizer of a quadratic function, and can then be rewritten by straightforward calculations as:
\begin{equation}
\label{eq:regression_coeff_proj}
  \beta^{n} = \E \left[ \phi(Z_{t_n}) \phi(Z_{t_n})^T  \right]^{-1} \E \left[  V_{\Delta t}(t_{n+1},Z_{t_{n+1}})  \phi(Z_{t_n}) \right],
\end{equation}
where we use the notation $ \phi = ( \phi_{1},\,\dots,\, \phi_{K})^T$, and where we assumed that $ \E \left[ \phi(Z_{t_n}) \phi(Z_{t_n})^T  \right]$ is invertible.

In order to estimate a solution to \eqref{eq:regression_coeff_proj} we rely on Monte Carlo simulations to approximate  expectations with finite sums.
Consider the training set $\{ \big( Z_{t_{n}}^m, Z_{t_{n+1}}^m\big)\}_{m=1}^{M}$ at time $t_{n}$ obtained by running $M \geq 1$ forward simulations of the process $Z$ from time $t_{0}=0$ to $t_{n+1}$.
$\beta^{n}$ defined in \eqref{eq:regression_coeff_proj} can then be estimated by
\begin{equation}
\label{eq:uncontrolled_coeff_sample}
\hat \beta^{n} \; = \left(\hat{\mathcal{A}}^M_n\right)^{-1} \frac{1}{M} \sum_{m=1}^M V_{\Delta t}(t_{n+1},Z^m_{t_{n+1}})\phi(Z^m_{t_n}),
\end{equation}
where we denote by $\hat{\mathcal{A}}^M_n$ the estimator $\frac{1}{M}\sum_{m=1}^M \phi(Z^m_{t_n})\phi(Z^m_{t_n})^T$ of the covariance matrix $\mathcal{A}_n=\E \left[\phi(Z_{t_n})\phi(Z_{t_n})^T \right] $.

The procedure presented above offers a convenient mean to approximate conditional expectations when the dynamics of the process $Z$ are uncontrolled.
When controlled, however, one has to account for the effect of the control on the conditional expectations either explicitly, via Control Randomization, or implicitly, via Regress-Later.

\vspace{3mm}
\noindent \textbf{Control Randomization}\\
In order to explicitly account for the effect of the control, one could directly introduce dependence on the control in the basis function.
This basic idea of Control Randomization consists in replacing in the dynamics of $Z$ the endogenous control by an exogenous control $(I_{t_n})_{0 \leq n \leq N}$, as introduced in \cite{Kharroubi2013}.
Its trajectories can then be simulated from time $t_0$ to time $t_N$.
Consider the training set $\{ Z_{t_n}^{m},{I}_{t_n}^m\}_{n=0,m=1}^{N,M}$, with $M\geq 1$, where ${I}_{t_n}^m$ are i.i.d. samples from a ``training distribution'' $\mu_n$ with support in $A$. The training set will be used to estimate the optimal $\beta^n$ coefficients for $n=0,\ldots,N-1$.
In the case of value iteration, $\{V_{\Delta t}^{}(t_{n+1},Z_{t_{n+1}}^{m}) \}_{m=1}^M$ is regressed against basis functions (which are, in this context, functions of the state and the control) evaluated at the points $\{ Z_{t_n}^{m},{I}_{t_n}^m \}_{m=1}^M$, as follows:
\[
	\E_{n,z}^a\big[V_{\Delta t}(t_{n+1},Z_{t_{n+1}})\big]\approx\sum_{k=1}^K\hat{\beta}^n_k\phi_k(z,a),
\]
with
\begin{align}
\hat{\beta}^{n}&=\argmin_{\beta\in\mathbb{R}^K}\left\{\sum_{m=1}^M \Big[V_{\Delta t} \big(t_{n+1},Z^m_{t_{n+1}} \big)-\sum_{k=1}^K\beta_k\phi_k(Z^m_{t_n},{I}^m_{t_n})\Big]^2\right\}\\
&\approx {\big(\hat{\mathcal{A}}^M_n \big)}^{-1}\frac{1}{M} \sum_{m=1}^M \Big[V_{\Delta t} \big(t_{n+1},Z^m_{t_{n+1}} \big) \phi(Z^m_{t_n},{I}^m_{t_n})\Big],
\end{align}
and where $\phi =(\phi_1,\ldots,\phi_K)^T$ and
\begin{equation}
\label{eq:def-estimator-covMatrix}
	\hat{\mathcal{A}}^M_n = \frac{1}{M}\sum_{m=1}^M \phi(Z^m_{t_n},I^m_{t_n})\phi(Z^m_{t_n},I^m_{t_n})^T
\end{equation}
 is an estimator of the covariance matrix $\mathcal{A}_n=\E[\phi(Z_{t_n},I_{t_n})\phi(Z_{t_n},I_{t_n})^T] $.

Notice that the basis functions take state and action variables as input in the case of Control Randomization-based method, i.e., their domain is $\R^d \times A$.
Also, observe that the estimated conditional expectation highly depends on the choice of the randomization for the control\footnote{Basically, different randomized controls may drive the process $Z$ to very different locations, and the estimations will suffer from inaccuracy on the states that have been rarely visited.}.

 An optimal feedback control at time $t_n$ given $Z_{t_n}=z$  is approximated by the expression (see Subsection \ref{sec:optimalControlSearching} for more practical details on the computation of the $\argmax_{} $):
\begin{equation}\label{eqn:CR_opt}
\hat \alpha_{t_n}(z) = \argmax_{a \in A} \left\{ f(z,a) \Delta t+  \sum_{k=1}^K\hat \beta_k^n\phi_k(z,a) \right\}.
\end{equation}

The value function at time $t_n$ is then estimated using Control Randomization method and value iteration procedure as
\[
\widehat V_{\Delta t}^{\text{CR}}(t_n,z) = f\big(z,\hat \alpha_{t_n}(z)\big) \Delta t+\sum_{k=1}^K\hat \beta^{n}_k {\phi}_k \big(z, \hat \alpha_{t_n}(z)\big), \quad z \in \R^d.
\]

Notice that Control Randomization can be easily employed in a performance iteration procedure by computing controls \eqref{eqn:CR_opt}, keeping in mind that at each time $t_n$ we need to re-simulate new trajectories $\{\tilde Z_{t_k}^m\}_{k=n,m=1}^{N,M}$ iteratively from the initial condition $\tilde{Z}^m_{t_n}=z$, using the estimated optimal strategies $(\hat \alpha_{t_k})_{k=n+1}^{N-1}$ to compute the quantities $ \sum_{k=n}^{N-1} f \big(t_k,\tilde Z^m_{t_k},\hat \alpha_{t_k}(\tilde Z^m_{t_k}) \big)+g(\tilde Z^m_{t_N}), \; \text{ for }  {1 \leq m \leq M}$.

\vspace{3mm}
\noindent {\textbf{Regress-Later}}\\
We present now a regress-later idea in which conditional expectation with respect to $Z_{t_n}$ is computed in two stages.
First, a conditional expectation with respect to $Z_{t_{n+1}}$ is approximated in a regression step by a linear combination of basis functions of $Z_{t_{n+1}}$.
Then, analytical formulas are applied to condition this linear combination of functions of future values on the present value $Z_{t_n}$.
For further details see \cite{Glasserman2002},  \cite{nadetal17} or \cite{balpal17}.
With this approach, the effect of the control is factored in implicitly, through its effect on the (conditional) distribution of $Z_{t_{n+1}}$ conditioned on $Z_{t_n}$.

Unlike the traditional Regress-Now method for approximating conditional expectations (which we discussed so far in the uncontrolled case and in Control Randomization), the Regress-Later approach, as studied in \cite{balpal17}, imposes conditions on basis functions:
\begin{assumption}
\label{laterass}
For each basis function $\phi_k$, $k=1, \ldots, K$, the conditional expectation
\[
\hat{\phi}_k^n(z,a) = \mathbb{E}^a_{n,z}[\phi_k(Z_{t_{n+1}})]
\]
can be computed analytically.
\end{assumption}


Using the Regress-Later approximation of the conditional expectation and recalling Assumption \ref{laterass} we obtain the optimal control $\alpha_{t_n}^m$ corresponding to the point $Z_{t_n}^m$, sampled independently from a ``training distribution'' $\mu_n$ (see Subsection \ref{sec:choiceGamma_n} for further details):
\[
\alpha_{t_n}^m=\argmax_{a\in A}\Bigl\{ f(Z_{t_n}^m,a) \Delta t+\sum_{k=1}^K\hat \beta^{n+1}_k  \hat \phi_k^n \big(Z_{t_n}^m,a\big)\Bigr\}.
\]
Notice that we are able to exploit the linearity of conditional expectations because

\begin{equation}
\label{eq:coeff_RL_sample}
\hat{\beta}^{n+1}=\argmin_{\beta\in\mathbb{R}^K}\left\{\sum_{m=1}^M\Big[V_{\Delta t}(t_{n+1},Z^m_{t_{n+1}})-\sum_{k=1}^K\beta_k\phi_k(Z^m_{t_{n+1}})\Big]^2\right\}
\end{equation}
is a constant once the training sets at times $t_k$, $k=n+1,\ldots,N,$ are fixed.\\

\noindent The value function at time $t_n$, is then estimated using Regress-Later method and value iteration procedure as
\[
\widehat V_{\Delta t}^{\text{RL}}(t_n,Z_{t_n}^m) = f(Z_{t_n}^m,\alpha_{t_n}^m) \Delta t+\sum_{k=1}^K\hat \beta^{n+1}_k \hat{\phi}_k^n \big(Z_{t_n}^m,\alpha_{t_n}^m\big).
\]
Notice that contrary to Control Randomization, Regress-Later does not require the training points to be distributed as $Z_{t_{n+1}}$ conditioned on $Z_{t_n}$ because the projection \eqref{eq:coeff_RL_sample} is an approximation to an expectation conditional to the measure $\mu_n$ which can be chosen freely to optimize the precision of the sample estimation.
On the other hand Regress-Later, similarly to Control Randomization, can be easily employed in a performance iteration procedure by generating forward trajectories at each time step.

\begin{remark}
Recall that the Regress-Later method uses training points that are i.i.d at each time step and independent across time steps.
Contrary to other Regression Monte Carlo approaches, Regress-Later does not require to use the information about the conditional distribution during the regression step as that is accounted for in the second step of the method, when conditional expectations are computed analytically.
\end{remark}

\subsubsection{Quantization}

We propose in this section a quantization-based algorithm to numerically solve control problems.
We may also refer to the latter as the Q-algorithm or Q in all the numerical examples considered in Section \ref{sec:numres}, where $Q$ stands for Quantization.
Let us first introduce some ingredients of Quantization, and then propose different ways of using them to approximate conditional expectations.

Let $(E,|.|)$ be a Euclidean space. We denote by $\hat\eps$ a $L$-quantizer of  an $E$-valued random variable $\eps$,
that is a discrete random variable on a grid $\Gamma$ $=$ $\{e_1,\ldots,e_L\}$ $\subset$ $E^L$ defined by
\beqs
\hat\eps &=& {\rm Proj}_\Gamma(\eps) \; = \; \sum_{\ell=1}^L e_\ell 1_{\eps \in C_i(\Gamma)},
\enqs
where $C_1(\Gamma)$, $\ldots$, $C_L(\Gamma)$ are the Voronoi cells corresponding to $\Gamma$, i.e., they form a Borel partition of $E$ satisfying
\beqs
C_\ell (\Gamma) & \subset & \Big\{ e \in E: |e-e_\ell | \; = \; \min_{j =1,\ldots,L} |e- e_j  | \Big\},
\enqs
and where ${\rm Proj}_\Gamma$ stands for the Euclidian projection on $\Gamma$.\\
The discrete law of $\hat\eps$ is then characterized by
\beqs
 p_\ell & =& \P[ \hat\eps = e_\ell ] \; = \; \P[ \eps \in C_\ell (\Gamma) ], \;\;\; \ell=1,\ldots,L.
\enqs
The grid of points $(e_\ell )_{\ell=1}^L$ which minimizes the $L^2$-quantization error $\| \eps - \hat\eps\|_{_2}$ leads to the so-called optimal $L$-quantizer, and can be obtained by a stochastic gradient descent method, known as Kohonen algorithm or competitive learning vector quantization (CLVQ) algorithm, which also provides  as a byproduct an estimation of the discrete law $(p_\ell )_{\ell=1}^L$.
We refer to \cite{pagetal04} for a description of the algorithm, and mention that for the normal distribution, the optimal grids and the weights of the Voronoi tesselations are precomputed for dimension up to 10 and are available on the website \url{http://www.quantize.maths-fi.com}.

In practice, optimal grids of the Gaussian random variable $\mathcal{N}_1(0,1)$ of dimension 1 with 25 to 50 points, have been used to solve the control problems considered in Section \ref{sec:numres}.

 We now propose two ways to approximate conditional expectations. The first approximation belongs to the family of the constant piecewise approximation, and the other one is an improvement on the first one, where the continuity of the approximation w.r.t. the control variable is preserved.

In the sequel, assume that for $n=0,\ldots,N-1$ we have a set $\Gamma_n$ of points in $\R^d$ that should be thought of as a training set used for estimating $V(t_n,\cdot)$.
See Subsection \ref{sec:choiceGamma_n} for more details on how to build $\Gamma_n$.

\vspace{3mm}
\noindent {\textbf{Piecewise constant interpolation}}

We assume here that we already have an estimate of $V_{\Delta t}(t_{n+1},\cdot)$, the value function at time $t_{n+1}$, for $n \in \{0,\ldots,N-1\}$, and we denote by $\widehat V^Q_{\Delta t}(t_{n+1},\cdot)$ the estimate.\\
The conditional expectation is then approximated as
\begin{align}
\label{approxEspCond1}
\mathbb{E}_{n,z}^a\big[\widehat{V}^Q_{\Delta t}\big(t_{n+1}, Z_{t_{n+1}}\big)\big] &\approx  \sum_{\ell =1}^L p_\ell \widehat{V}^Q_{\Delta t}\Big(t_{n+1}, \text{Proj}_{\Gamma_{n+1}} \big(G_{\Delta t}(z,a,e_\ell)\big) \Big), \quad \text{ for } z \in \Gamma_{n},
\end{align}
where:
\begin{itemize}
  \item $G_{\Delta t}$ is defined, using the notations introduced in \eqref{eq:euler}, as
  \begin{equation}
  \label{def:G_deltat}
    G_{\Delta t}(z,a,\eps)=z+ b(z, a) \Delta t + \sigma_0(z,a ) \sqrt{\Delta t} \; \eps.
  \end{equation}
  \item $\text{Proj}_{\Gamma_{n}}(.)$ stands for the Euclidean projection on $\Gamma_{n}$.
  \item $\Gamma= \{ e_1, \ldots, e_L \big\}$ and $\big\{ p_\ell \big\}_{1 \leq \ell \leq L}$ are respectively the optimal $L$-quantizer and its associated sequence of weights of the exogenous noise $\eps$. See above for more details.
\end{itemize}

 An optimal feedback control at time $t_n$ and point $z \in \Gamma_n$  is approximated by the expression (see Subsection \ref{sec:optimalControlSearching} for more practical details on the computation of the $\argmax_{} $):
\begin{equation}\label{eqn:Q_Contopt}
\hat \alpha^Q_{t_n}(z) = \argmax_{a \in A} \left\{ f(z,a) \Delta t+   \sum_{\ell =1}^L p_\ell \widehat{V}^Q_{\Delta t}\Big(t_{n+1}, \text{Proj}_{\Gamma_{n+1}} \big(G_{\Delta t}(z,a,e_\ell)\big) \Big) \right\}.
\end{equation}

The value function at time $t_n$, is then estimated using the piecewise constant approximation and value iteration procedure as
\[
\widehat V_{\Delta t}^{\text{Q}}(t_n,z) = f\big(Z_{t_n}^m,\hat \alpha^Q_{t_n}(z) \big) \Delta t+  \sum_{\ell =1}^L p_\ell \widehat{V}^Q_{\Delta t}\bigg(t_{n+1}, \text{Proj}_{\Gamma_{n+1}} \Big(G_{\Delta t}\big(z, \hat \alpha^Q_{t_n}(z),e_\ell \big )\Big) \bigg).
\]

\begin{remark}
Clearly, the constant piecewise approximation can be easily extended to control problems of all dimensions $d\geq 1$.
However the latter is, in most cases, not continuous w.r.t. the control variable since it remains constant on each Voronoi cells (see, e.g., Figure \ref{fig:DiscApprox} p.\pageref{fig:DiscApprox}).
As a result, the optimization process over the control space suffers from high instability and inaccuracy, which implies a poor estimation of the value function $V(t_n,\cdot)$.
\end{remark}

\vspace{3mm}
 \noindent {\textbf{Semi-linear interpolation}}
Once again, we assume here that we already have $\widehat V^{\text{Q}}_{\Delta t}(t_{n+1},\cdot)$, an estimate of the value function at time $t_{n+1}$, with $n \in \{0,\ldots,N-1\}$, and wish to provide an estimation of the conditional expectation in the particular case where the controlled process lies in dimension $d$=1.
Consider the following piecewise linear approximation of the conditional expectation, which is continuous w.r.t. the control variable $a$:
\begin{align}
\label{approxEspCond2}
\mathbb{E}_{n,z}^a\big[\widehat{V}^\text{Q}_{\Delta t}\big(t_{n+1}, Z_{t_{n+1}} \big)\big] &\approx \sum_{ \ell =1}^L p_\ell \Big[\lambda^{e_\ell,z}_a\widehat{V}^\text{Q}_{\Delta t}\big(t_{n+1},z_+\big) +(1-\lambda^{e_\ell,z}_a)\widehat{V}^\text{Q}_{\Delta t}\big(t_{n+1}, z_- \big) \Big], \quad \text{ for } z \in \Gamma_n,
\end{align}
where for all $\ell=1,\ldots,L$, $z_-$ and $z_+$ are defined as follows:
\begin{itemize}
  \item $z_-$ and $z_+$ are the two closest states in $\Gamma_{n+1}$ from $ G_{\Delta t}(z,a,e_\ell) $, such that $z_- < G_{\Delta t}(z,a,e_\ell) < z_+$, if such states exist; and, in this case, we define $\lambda^{e_\ell,z}_a = \frac{G_{\Delta t}(z,a,e_\ell)-z_-}{z_+-z_-}$.
  \item Otherwise, $z_-$ and $z_+$ are equal and defined as the closest state in $\Gamma^Z_{n+1}$ from $ G_{\Delta t}(z,a,e_\ell) $  and we define $\lambda^{e_\ell,z}_a=1$.
\end{itemize}

\begin{remark}
  This second approximation is continuous w.r.t. the control variable, which brings stability and accuracy to the optimal control task (see Subsection \ref{sec:optimalControlSearching}), and also ensures an accurate estimate of the value function at time $t_n$.
We will mainly use this approximation in the numerical tests (see Section \ref{sec:numres}).
\end{remark}

\begin{remark}
  Although the dimension $d=1$ plays a central role to define clearly the states $z_-$ and $z_+$ in \eqref{approxEspCond2}, the semi-linear approximation can actually be generalized to a certain class of control problems of dimension greater than 1, using multi-dimensional Quantization (see, e.g., the comments on the Q-algorithm designed to solve the Portfolio Optimization example, in Subsection \ref{rem:QalgoPortfolioLiq}).
However, it is not well-suited to solve numerically general control problems in dimension greater than 1.
For these cases, other interpolating methods such as the use of Gaussian processes are more appropriated (see, e.g., \cite{lud18} for an introduction on the use of  Gaussian processes in Regression Monte Carlo).
\end{remark}

 The optimal feedback control at time $t_n$ and point $z \in \Gamma_n$  is approximated as (see Subsection \ref{sec:optimalControlSearching} for more practical details on the computation of the $\argmax_{}$):
\begin{equation}\label{eqn:CR_opt2}
\hat \alpha^\text{Q}_{t_n}(z) = \argmax_{a \in A} \left\{ f(z,a) \Delta t+   \sum_{\ell =1}^L p_\ell \Big[\lambda^{e_\ell,z}_a\widehat{V}^\text{Q}_{\Delta t}\big(t_{n+1},z_+\big) +(1-\lambda^{e_\ell,z}_a)\widehat{V}^\text{Q}_{\Delta t}\big(t_{n+1}, z_- \big) \Big] \right\}.
\end{equation}

The value function at time $t_n$ is then estimated using the semi-linear approximation and value iteration procedure as
\[
\widehat V_{\Delta t}^{\text{Q}}(t_n,z) = f\big(z,\hat \alpha^\text{Q}_{t_n}(z) \big) \Delta t+  \sum_{\ell=1}^L p_\ell \Big[\lambda^{e_\ell,z}_{\hat \alpha^\text{Q}_{t_n}(z)}\widehat{V}^\text{Q}_{\Delta t}\big(t_{n+1},z_+\big) +(1-\lambda^{e_\ell,z}_{\hat \alpha^\text{Q}_{t_n}(z)})\widehat{V}^\text{Q}_{\Delta t}\big(t_{n+1}, z_- \big) \Big],
\]
where $z_+$ and $z_-$ are defined using the control $\hat \alpha^\text{Q}_{t_n}(z)$. See \eqref{approxEspCond2} for their definitions.

\subsection{Training points design}
\label{sec:choiceGamma_n}

 We discuss here the choice of the training measure $\mu$ and the sets $(\Gamma_{n})_{n=0,\ldots,N-1}$ used to compute the numerical approximations in Regression Monte Carlo and Quantization.
Two cases are considered in this section.
The first one is a knowledge-based selection, relevant when the controller knows with a certain degree of confidence where the process has to be driven in order to optimize her reward functional.
The second case, on the other hand, is when the controller has no idea where or how to drive the process to optimize the reward functional.

 \subsubsection{Exploitation only strategy} In the knowledge-based setting there is no need for exhaustive and expensive (in time mainly) exploration of the state space, and the controller can directly choose training sets $\Gamma$ constructed from distributions $\mu$ that assign more points to the parts of the state space where the optimal process is likely to be driven.

 In practice, at time $t_n$, assuming we know that the optimal process is likely to stay in the ball centered around the point $m_n$ and with radius $r_n$, we chose a training measure $\mu_n$ centered around $m_n$ as, for example $\mathcal{N}(m_n,r_n^2)$, and build the training set as sample of the latter.
In the Regress-Later setting this can be done straightforwardly, while Control Randomization requires one to select a measure for the random control such that the controlled process $Z$ is driven in such area of the state space.

Taking samples according to $\mu$ to build grids makes them random.
Another choice, which we used in the Quantization-based algorithm, is to use the (deterministic) optimal grid of $\mathcal{N}(m_n,\sigma_n^2)$ with reduced size (typically take 50 points for a problem in dimension 1, 250 for one of dimension 2 when $\sigma_n^2=1$,\ldots), which can be found at \url{www.quantize.maths-fi.com}, to reduce the size of the training set and alleviate the complexity of the algorithms.

\begin{remark}
  As the reader will see, we chose the training sets based on the ``exploitation only strategy'' procedure, i.e. by guessing where to drive optimally the process, when solving the Liquidation Problem introduced  in Subsection \ref{sec:portfolioLiquidation}.
\end{remark}

 \subsubsection{Explore first, exploit later}
 \label{sec:explorefirstregresslater}

 \noindent {\textbf{Explore first:}}
If the agent has no idea of where to drive the process to receive large rewards, she can always proceed to an exploration step to discover favorable subsets of the state space.
To do so, the $\Gamma_{n}$, for $n=0,\ldots,N-1$, can be built as uniform grids that cover a large part of the state space, or $\mu$ can be chosen uniform on such domain.
It is essential to explore far enough to have a well understanding of where to drive and where not to drive the process.

\vspace{3mm}
\noindent {\textbf{Exploit later:}}
The estimates for the optimal controls at time $t_n,\; n=0,\ldots, N-1$, that come up from the \textit{Explore first} step, are relatively good in the way that they manage to avoid the wrong areas of state space when driving the process.
However, the training sets that have been used to compute the estimated optimal control are too sparse to ensure accuracy on the estimation.
In order to improve the accuracy, the natural idea is to build new training sets by simulating $M$ times the process using the estimates on the optimal strategy computed from the \textit{Explore first} step, and then proceed to another estimation of the optimal strategies using the new training sets.
This trick can be seen as a two steps algorithm that improves the estimate of the optimal control.

\begin{remark}
In Control Randomization, multiple runs of the method are often needed to obtain precise estimates, because the initial choice of the dummy control could drive the training points far from where the optimal control would have driven them.
In practice, after having computed an approximate policy backward in time, such policy is used to drive $M$ simulations of the process forward in time, which in turn produce control paths that can be fed as a random controls in a new backward procedure, leading to more accurate results.
\end{remark}

\begin{remark}
We applied the ``explore first, exploit later'' idea to solve the Portfolio Optimization problem introduced in Subsection \ref{sec:portfolioLiquidation}.
\end{remark}

\subsection{Optimal control searching}
\label{sec:optimalControlSearching}
Assume in this section that we already have the estimates $\widehat V_{\Delta t}(t_k,\cdot)$ for the value function at time $t_k$, for $k=n+1,\ldots,N$, and want to estimate $V(t_n,\cdot)$ the value function at time $t_n$.

The optimal control searching task consists in optimizing the function\footnote{often referred to as the $Q$-function, or action-value function, in the reinforcement learning literature.
Be aware that $Q$ stands here for the "Quality" of an action taken in a given state, and in particular does not refer to Quantization.}
\[
\widehat{Q}_n:(z,\cdot) \mapsto f(z,a) \Delta t + \hat{\E}_{n,z}^a\big[\widehat V_{\Delta t}(t_{n+1},Z_{t_{n+1}})\big]
\]
over the control space $A$, for each $z \in \Gamma_{n}$, and where we denote by $\hat{\E}_{n,z}^a\big[\widehat V_{\Delta t}(t_{n+1},Z_{t_{n+1}})\big]$ an approximation of $\E_{n,z}^a\big[\widehat V_{\Delta t}(t_{n+1},Z_{t_{n+1}})\big]$ using Regress-Later, or Control Randomization or Quantization-based methods (see Subsection \ref{sec:approxexpect}).
Once again, we remind that importance of this task is motivated by the dynamic programming principle stating that for all $n=0,\ldots,N-1$, we can approximate the value function at time $n$ as follows
\begin{equation}
\widehat V_{\Delta t}(t_n,z)= \Sup_{a \in A}  \; \widehat Q_n(z,a),
\end{equation}
where $\widehat{V}_{\Delta t}(t_n,\cdot)$ is our desired estimate of the value function at time $n$.

\subsubsection{Low cardinality control set}
In the case where the control space $A$ is discrete (with a relatively small cardinality), one can solve the optimization problem by an exhaustive search over all the available controls without compromising the computational speed.

\begin{remark}
Note that in the case where the control space is continuous, one can always discretize the latter in order to rely on the effectiveness of extensive search to solve the optimal control problem.
However, the control space discretization brings an error.
So the control might have to include a high number of points in the discretization in order to reduce the error thereby causing a considerable slow down of the computations.
\end{remark}

\subsubsection{High cardinality/continuous control space}

If we assume differentiability almost everywhere, as follows from the semi-linear approximation in Quantization, and most choices of basis functions in Regression Monte Carlo, we can carry on the optimization step by using some gradient-based algorithm for optimization of differentiable functions.
Actually, many optimizing algorithms (Brent, Golden-section Search, Newton gradient-descent,\ldots) are already implemented in standard libraries of most programming languages like Python (see, e.g., package \texttt{scipy.optimize}), Julia (see, e.g., package \texttt{Optim.jl}), C and C++ (see, e.g., package \texttt{NLopt}).

\begin{remark}
  \label{rq:instantaneousFinding}
When the control is of dimension 1, polynomials of order smaller than 5 are employed as basis functions in Regression Monte Carlo as well as for the running reward $f$.
 The optimal control can then be computed analytically as a function of the regression coefficients, since every polynomial equation of order smaller than 4 can be solved by radicals.
\end{remark}

 Concretely, in all the examples considered in Section \ref{sec:numres}, we used the Golden-section Search or the Brent methods when testing Quantization-based algorithm to find the optimal controls at each point of the grids.
These algorithms were very accurate to find the optimal controls, and we made use of Remark \ref{rq:instantaneousFinding} to find the optimal controls using the Regress-Later-based algorithm.

\subsection{Upper and lower bounds}

After completing the backward procedure, we can compute an unbiased estimation of the value of the control policy by using Monte Carlo simulations and sample average.
Assume already computed (or simply available) the matrix of regression coefficients, in the case of Regression Monte Carlo, and discrete probability law $\hat p$ for Quantization, we can use this information to implicitly compute the control and simulate forward many trajectories of the controlled process starting from a common initial condition.
We can then evaluate the average performance measure by computing the sample average of the rewards collected on each trajectory.
Denoting such approximation by $\doublehat V_{\Delta t}(0,z)$, and recalling that by definition $J_{\Delta t}(0,z,\alpha)\le J_{\Delta t}(0,z,\alpha^*)$, for all $ \alpha \in  \mathcal{A}_{\Delta t}$ and where $\alpha^*$ represents the optimal control process; it holds $\doublehat V_{\Delta t}(0,z) \leq V_{\Delta t}(0,z),$ for $z \in \R^d$.

The argument above implies that, neglecting the time-discretization error, we obtain a lower bound for $V_{\Delta t}(0,\cdot)$ by evaluating the estimated policy.
On the other hand, see \cite{PHL17}, based on \cite{rog02}, to get an upper bound of the value function via duality.

\subsection{Pseudo-codes}
\label{sec:pseudocodes}

In this section, we present the pseudo-code for the three approaches presented in the previous sections.
For simplicity, we will only show the algorithms designed using value iteration procedure.
However, the performance iteration update rule can be substituted in the codes below provided that forward simulations are run to obtain a pathwise realization of the controlled process and associated rewards.

\subsubsection{Pseudo-code for a Regress-Later-based algorithm}

 We present in Algorithm \ref{alg:RLMC} a pseudo-code to estimate $V_{\Delta t}(t_n,\cdot)$, for $n=0,\ldots,N-1$, using Value Iteration and based on Regress-Later method.
For $n=0,\ldots,N-1$, we denote by $\hat V^{\text{RL}}_{\Delta t} (t_n,\cdot)$ the derived estimation of $V_{\Delta t}(t_n,\cdot)$, and will refer to it as the RLMC algorithm in the numerical tests presented in Section \ref{sec:numres}.

 Note that we use the same training measure $\mu$ at each time step so that there is only one covariance matrix to estimate  (since $\mathcal{A}_{t_n}$ is the same for all $n=0,\dots,N-1$). Denote by $ \hat{\mathcal{A}}^M$ the estimator, as defined in~\eqref{eq:def-estimator-covMatrix}.
\begin{algorithm}
 \caption{Regress-Later Monte Carlo algorithm (RLMC) - Value iteration
 \label{alg:RLMC}}
 \hbox{\textbf{Inputs:}}
   \begin{itemize}
    \item $M$: number of training points,
    \item $\mu$: distribution of training points,
    \item $K$: number of basis functions,
    \item $\{\phi_k\}_{k=1}^K$: family of basis functions,
  \end{itemize}
 \begin{algorithmic}[1]
\State Estimate the covariance matrix $\hat{\mathcal{A}}^M$
 \State Generate i.i.d. training points $\{{Z}_{t_N}^{m}\}_{m=1}^{M}$ accordingly to the distribution $\mu$.
\State Initialize the value function $\hat V^{\text{RL}}_{\Delta t}(t_N, Z_{t_N}^m)  =g({Z}_{t_N}^{m}), \quad \forall m=1, \dots, M$
 \For{$n=N-1$ to $0$}
    \State  $\hat \beta^{n}=\hat{\mathcal{A}}^{-1}_M\frac{1}{M}\sum_{m=1}^{M}\Big[\hat V^\text{RL}_{\Delta t}(t_{n+1},{Z}_{t_{n+1}}^{m}) \phi({Z}_{t_{n+1}}^{m})\Big]$
     \State Generate a new layer of i.i.d. training points $\{{Z}_{t_n}^{m}\}_{m=1}^{M}$
           accordingly to the distribution $\mu$.
     \State For all $m=1, \dots, M$ do
    \[
      \hat V^{\text{RL}}_{\Delta t}(t_n, Z_{t_n}^m)
  =\sup_{a \in A} \Big\{f({Z}_{t_n}^{m},a) \Delta t+\sum_{k=1}^K \hat \beta^{n}_k \hat {\phi}_k^n({Z}_{t_n}^{m}, a)\Big\}
    \]
\EndFor
\State \textbf{Evaluate the policy to obtain $\doublehat V^{\text{\normalfont RL}}_{\Delta t}$}
\end{algorithmic}
\hbox{\textbf{Outputs:} $\{\hat{\beta}_{k}^n\}_{n,k=1}^{N,K}$, $\doublehat V^{\text{RL}}_{\Delta t}(0,z)$ for $z \in \R^d$.}
 \end{algorithm}

\subsubsection{Pseudo-code for a Control Randomization-based algorithm}

 We present in Algorithm \ref{alg:CR} a pseudo-code to estimate $V_{\Delta t}(t_n,\cdot,\cdot)$, for $n=0,\ldots,N-1$, using Value Iteration and based on Control Randomization method.
For $n=0,\ldots,N-1$, we denote by $\hat V^{\text{CR}}_{\Delta t} (t_n,\cdot)$ the derived estimation of $V_{\Delta t}(t_n,\cdot)$, and will refer to it as the CR algorithm in the numerical tests presented in Section \ref{sec:numres}.

 \begin{algorithm}
 \caption{Control Randomization algorithm (CR) - Value iteration
 \label{alg:CR}}
 \hbox{\textbf{Inputs:}}
   \begin{itemize}
    \item $M$: number of training points,
    \item $\mu$: initial distribution of dummy control,
    \item $K$: number of basis functions,
    \item $\{\phi_k\}_{k=1}^K$: family of basis functions,
  \end{itemize}
 \begin{algorithmic}[1]
\State Estimate the covariance matrix $\hat{\mathcal{A}}^M$
 \State Generate $m$ trajectories, $\{ Z_{t_n}^{m},I_{t_n}^{m}\}_{n=0,m=1}^{N,M}$, where $Z_{t_n}^{m}$ is driven by $I_{t_n}^{m}$, and the  $I_{t_n}^{m}$ are i.i.d with distribution $\mu$.
\State Initialize the value function $\hat  V^{\text{CR}}_{\Delta t}(t_N,   Z_{t_N}^m)  = g( {Z}_{t_N}^{m}), \quad m=1, \dots, M$
  \For{$n=N-1$ to $0$}
    \State  $\hat \beta^{n}=(\hat{\mathcal{A}}^M)^{-1}\frac{1}{M}\sum_{m=1}^{M}\Big[\hat V^{\text{CR}}_{\Delta t}(t_{n+1},{Z}_{t_{n+1}}^{m}) \phi({Z}_{t_{n}}^{m},I_{t_n}^{m})\Big]$
     \State For all $m=1, \dots, M$ do
    \[
      \hat V^{\text{CR}}_{\Delta t}(t_n, Z_{t_n}^m)
  =\sup_{a \in A} \Big\{f({Z}_{t_n}^{m},a) \Delta t+\sum_{k=1}^K \hat\beta^{n}_k\phi_k({Z}_{t_n}^{m}, a)\Big\}
    \]
\EndFor
\State \textbf{Evaluate the policy to obtain $\doublehat V^{\text{\normalfont CR}}_{\Delta t}$}
\end{algorithmic}
\hbox{\textbf{Outputs:} $\{\hat{\beta}_{k}^n\}_{n=0,k=1}^{N,K}$, $\doublehat V^{\text{CR}}_{\Delta t}(0,z)$ for $z \in \R^d$.}
 \end{algorithm}

\subsubsection{Pseudo-code for a Quantization-based algorithm}
 We present in Algorithm \ref{alg:QUANT} a pseudo-code to estimate $V_{\Delta t}(t_n,\cdot,\cdot)$, for $n=0,\ldots,N-1$, using value iteration procedure and based on Quantization method.
For $n=0,\ldots,N-1$, we denote by $\hat V^{\text{Q}}_{\Delta t} (t_n,\cdot)$ the derived estimation of $V_{\Delta t}(t_n,\cdot)$, and will refer to it as the Q-algorithm in the numerical tests presented in Section \ref{sec:numres}.

Note that we made use of a piecewise constant approximation of conditional expectations to approximate $\hat V^{\text{Q}}_{\Delta t}(t_n,\cdot)$ in order to keep the algorithm simple.
Also, note that, as said previously, in most of the numerical tests run in Section \ref{sec:numres}, we will use optimal grids available at \url{www.quantize.maths-fi.com} and will take $L=25$ to $50$ points for the size of the optimal grid of the Gaussian noise $\eps$.

\begin{algorithm}
  \caption{Quantization algorithm (Q) - Value iteration
    \label{alg:QUANT}}
  \hbox{\textbf{Inputs:}}
  \begin{itemize}
    \item $\Gamma_k$, $k=0,\ldots,N$: grids of training points in $\R^d$,
    \item $\Gamma= \{e_1, \ldots,e_L\}$\;, $(p_\ell)_{1\leq \ell \leq L}$ : the L-optimal grid of the exogenous noise $\eps$, and its associated weights.
  \end{itemize}
  \begin{algorithmic}[1]
    \State Initialize the estimated value function at time $N$: $\hat V^{\text{Q}}_{\Delta t}(t_N, z)  = g(z), \quad \forall z \in \Gamma_{N}$.
    \For{$n=N-1$ to $0$}
    \State  Estimate the value function at time $t_n$ as follows:
    \begin{equation}
    \label{progdynQ}
    \widehat V^\text{Q}_{\Delta t}(t_{n}, z)=\Max_{a \in A}\Bigg[  f(z,a) \, \Delta t +\sum_{\ell=1}^L p_\ell \widehat{V}^\text{Q}_{\Delta t}\Big(t_{n+1}, \text{Proj}_{\Gamma_{n+1}} \big(G_{\Delta t}(z,a,e_\ell )\big) \Big) \Bigg], \quad \forall z \in \Gamma_{n}.
    \end{equation}
    \State Compute the optimal strategies $\hat \alpha(t_n,z), \; {z \in \Gamma_n}$, as maximizer of \eqref{progdynQ}:
    \begin{equation*}
    \hat \alpha(t_n,z) \in \argmax_{a \in A}\bigg[  f(z,a) \, \Delta t + \sum_{\ell=1}^L p_\ell \widehat{V}^{\text{Q}}_{\Delta t}\Big(t_{n+1}, \text{Proj}_{\Gamma_{n+1}} \big(G_{\Delta t}(z,a,e_\ell)\big) \Big) \bigg], \quad \forall z \in \Gamma_{n}.
    \end{equation*}
    \EndFor
    \State \textbf{Evaluate the policy to obtain $\doublehat V^{\text{\normalfont Q}}_{\Delta t}$}
  \end{algorithmic}
  \hbox{\textbf{Outputs:}  $\big( \hat \alpha(t_n, z)\big) _{z \in \Gamma_n, 0\leq n \leq N-1} $ ,  $\big( \doublehat V^{\text{Q}}_{\Delta t}(0,z) \big)_{z \in \Gamma_0}$. }
\end{algorithm}

 Table
\section{Applications and numerical results}
\label{sec:numres}

\subsection{Portfolio Optimization under drift uncertainty}
\label{sec:portfolioLiquidation}
\subsubsection{The model}

We consider a financial market model with one risk-free asset, assumed to be equal to one, and $d$  risky assets of price process $S$ $=$ $(S^1,\ldots,S^d)$ governed
\begin{equation*}
  dS_t = {\rm diag}(S_t)( \beta_t dt + \sigma dB_t^0) , \quad S_0=s_0 \in \R^d,
\end{equation*}

where $B^0$ is a $d$-dimensional Brownian motion on a filtered probability space $(\Omega,\Fc,\F,\P^0)$, $\sigma$  is the $d\times d$ invertible  matrix volatility coefficient, assumed to be known and constant.
However,  the drift $(\beta_t)$ of the asset (which is typically a diffusion process governed by another independent Brownian motion $B$) is unknown and  unobservable like the Brownian motion $B^0$.
The agent can actually only observe the stock prices $S$, and we denote by  $\F^S$ the filtration generated by the price process $S$, which should be view as the available information.

In this context, we shall consider two important classes of optimization problems in finance:
\begin{itemize}
\item[(1)] {\it Portfolio Liquidation.}
We consider the problem of an agent (trader) who has to liquidate a large number  $y_0$ of shares in some asset (we consider one stock, $d=1$) within a finite time $T$, and faces  execution costs and market price impact.
In contrast with frictionless Merton problem, we do not consider mark-to-market value of the portfolio and instead consider  separately the amount on the cash account and the inventory $Y$, i.e., the position or number of shares held at any time.
The strategy of the agent is then described by a real-valued  $\F^S$-adapted process $\alpha$, representing the velocity  at which she buys ($\alpha_t$ $>$ $0$) or sells ($\alpha_t$ $<$ $0$) the asset, and the inventory is thus given by
\beqs
Y_t &=& y_0 + \int_0^t \alpha_u du, \;\;\; 0\leq t \leq T.
\enqs
The objective of the trader is to minimize over $\alpha$ the total liquidation cost
\beqs
J_1(\alpha) &=&  \E^0 \big[ \int_0^T \alpha_t \big( S_t +   f(\alpha_t) \big) dt +   \ell(Y_T)  \big]
\enqs
where $f(.)$ is an increasing function with $f(0) = 0$,  representing a temporary price impact, and $\ell(.)$ is a loss function, i.e., a convex function with $\ell(0)$ $=$ $0$, penalizing the trader when she does not succeed to liquidate all her shares.

\item[(2)] {\it Portfolio Selection.}  The set $\Ac$ of portfolio strategies, representing the amount  invested in the assets,  consists in all $\F^S$-adapted processes $\alpha$ valued in some set $A$ of $\R^d$, and satisfying $\int_0^T |\alpha_t|^2 dt$ $<$ $\infty$.
The dynamics of wealth process $X$ $=$ $X^\alpha$ associated to a portfolio strategy $\alpha$ is then governed by
\beqs
dX_t &=& \alpha_t S_t^{-1}dS_t \\
 &= &   \alpha_t . \beta_t dt + \alpha_t\trans\sigma dB_t^0, \;\;\; X_0 \; = \; x_0 \in \R,
\enqs
and as in Merton Portfolio Selection problem, the objective of the agent is to maximize over portfolio strategies the utility of terminal wealth
\beqs
J_2(\alpha) &=& \E^0[ U(X_T)],
\enqs
where $U$ is a utility function on $\R$, e.g., CARA function $U(x)$ $=$ $-\exp(-p x)$, $p$ $>$ $0$.
\end{itemize}

Let us show how one can reformulate the above problems into a McKean-Vlasov type problem under partial observation and common noise as described in Section \ref{intro}.
We first introduce the so-called probability reference $\P$, which makes the observation price process a martingale.
Let us then define the process
\beqs
Z_t &=& \exp \big( - \int_0^t \sigma^{-1}\beta_u dB_u^0  - \frac{1}{2} \int_0^t |\sigma^{-1}\beta_u|^2 du \big), \;\;\; 0 \leq t \leq T,
\enqs
which is a $(\P^0,\F)$-martingale (under suitable integrability conditions on $\beta$), and defines a probability measure $\P$ $\sim$ $\P^0$ through  its density:
$\frac{d\P}{d\P^0}\Big|_{\Fc_t}$ $=$ $Z_t$, and under which the process
\beqs
W_t^0 &=& B_t^0 + \int_0^t \sigma^{-1} \beta_u du, \;\;\;  0 \leq t \leq T,
\enqs
 is a $(\P,\F)$-Brownian motion by Girsanov's theorem, and the dynamics of $S$ is
 \beqs
dS_t &=& {\rm diag}(S_t) \sigma dW_t^0.
 \enqs
Notice that $\F^S$ $=$ $\F^0$ the filtration generated by $W^0$.
We also denote by $L_t$ $=$ $1/Z_t$ the $(\P,\F)$-martingale governed by
\beqs
dL_t &=& L_t  \sigma^{-1} \beta_t . dW_t^0.
\enqs

Next, we use Bayes formula and rewrite the gain (resp. cost) functionals of our two portfolio optimization problems as
\beqs
J_1(\alpha) &=&  \E \big[ \int_0^T L_t \alpha_t ( S_t + f(\alpha_t)) dt +   L_T  \ell(Y_T )  \big] \\
&=& \E \big[ \int_0^T \bar L_t^0 \alpha_t ( S_t +  f(\alpha_t)) dt +   \bar L_T^0 \ell(Y_T)  \big] \; \\
&=& \;
\E \big[ \int_0^T \bar L_t^0 \alpha_t ( \bar S_t^0 + f(\alpha_t) ) dt +   \bar L_T^0 \ell(\bar Y_T^0)  \big]\\
J_2(\alpha) & = & \E \big[ L_T U(X_T) \big]  \; = \; \E \big[ \bar L_T^0 U(X_T) \big] \; = \; \E \big[ \bar L_T^0 U(\bar X_T^0) \big]
\enqs
where  $\bar L_t^0$ $=$ $\E[L_t |W^0]$ $=$ $\int \ell \P_{L_t}^{W^0}(d\ell)$, $\bar X_t^0$ $=$ $\E[X_t |W^0]$ $=$ $\int x \P_{X_t}^{W^0}(dx)$ $=$ $X_t$,
$\bar Y_t^0$ $=$ $\E[Y_t |W^0]$ $=$ $\int y \P_{Y_t}^{W^0}(dy)$ $=$ $Y_t$
$\bar S_t^0$ $=$ $\E[S_t |W^0]$ $=$ $\int s \P_{S_t}^{W^0}(ds)$ $=$ $S_t$,  and we used the law of conditional expectations and  the fact that $S$, $X$ and $Y$ are $\F^0$-adapted.
This formulation  of the functional  $J_1$ (resp. $J_2$) fits into the MKV framework of Section \ref{intro} with state variables $(X,L,\beta)$ (resp. $(Y,S,L,\beta)$)

We now consider the particular case when $\beta$  is an $\Fc_0$-measurable  random variable distributed according to some probability distribution $\nu(db)$: this corresponds to a Bayesian point of view when the agent's belief  about the drift  is modeled by a prior distribution.
In this case, let us show how our partial observation  problem can be embedded into a finite-dimensional full observation Markov control problem.
Indeed, by noting that $\beta$ is independent of the Brownian motion $W^0$ under $\P$, we have
\beqs
\bar L_t^0 &=& \E \big[ \exp \big( \sigma^{-1}\beta . W_t^0 - \frac{1}{2} |\sigma^{-1} \beta|^2 t \big) \big| W^0 \big]  \; = \; F(t,W_t^0),
\enqs
where
\beqs
F(t,w) & =& \int \exp\big( \sigma^{-1}b . w - \frac{1}{2}|\sigma^{-1}b|^2 t  \big) \nu(db).
\enqs
Hence, the functionals $J_1$ and $J_2$ can  be written as
\beq
J_1(\alpha) &=& \E \big[ \int_0^T F(t,W_t^0)  \alpha_t ( S_t + f( \alpha_t)) dt +  F(T,W_T^0) \ell(Y_T)  \big] \\
J_2(\alpha) & = & \E \big[ F(T,W_T^0) U(X_T) \big] \label{J1Markov}.
\enq
We are then reduced to a $(\P,\F^0)$-control problem with state variables $(W^0,X)$ for  problem (1) and $(W^0,S,Y)$ for problem (2) with dynamics under $\P$:
\beq
\label{model:continousT}
dS_t &=& {\rm diag}(S_t) \sigma dW_t^0, \quad S_0=s_0 \in (\R_+)^d \\
dX_t &=&  \alpha_t\trans\sigma dW_t^0, \quad X_0=0  \\
dY_t &=& \alpha_t dt, \quad Y_0=y_0 \in \R_+.
\enq

\begin{remark}
{\rm
Another example of partial observation for the drift $\beta$ is the case when it is modeled by  a linear Gaussian process.
This would lead to  the well-known Kalman-Bucy filter, hence to a finite-dimensional control problem.
However, for general unobserved drift process $\beta$,  we fall into an infinite dimensional control problem involving the filter process.
}
 \end{remark}

\subsubsection{Numerical results}
\label{numericalresults}

Let us now illustrate numerically the impact of uncertain Bayesian drift on the Portfolio Liquidation problem and the Portfolio Selection problem in dimension $d=1$, by considering a Gaussian prior distribution $\beta$ $\sim$ $\nu$ $=$  $\Nc( b_0, \gamma_0^2)$, with $b_0 \in \R$ and $\gamma_0 >0$.
In this case, $F$ is explicitly given by:
  \begin{equation*}
F(t,w)=\frac{\sigma}{ \sqrt{\sigma^2+\gamma_0^2 t}}\exp\Big(\frac{1}{2(\sigma^2+\gamma_0^2t)}(-b_0^2t+2b_0\sigma w+\gamma_0^2w^2) \Big).
\end{equation*}

{\bf 1. Portfolio Liquidation}. Let us first  consider the Portfolio Liquidation problem (1) with a linear price impact function $f(a)$ $=$ $\gamma a$, $\gamma$ $>$ $0$, and a quadratic loss function
$\ell(y)$ $=$ $\eta y^2$, $\eta$ $>$ $0$.
The optimal trading rate is given by (see  \cite{pha16})
\beqs
\alpha_t^* &=& - \frac{ Y_t^*}{T-t + \gamma/\eta} + \frac{1}{2\gamma} \Big(  \frac{ 1}{T-t + \gamma/\eta} \int_t^T \E^0[ S_u | \Fc^S_t ] du  - S_t   \Big)
\enqs
where $Y^*$ is the associated inventory with feedback control $\alpha^*$: $dY_t^*$ $=$ $\alpha_t^* dt$, $Y_0^*$ $=$ $y_0$.
Since we consider a Gaussian prior $\Nc( b_0, \gamma_0^2)$ for $\beta$, the optimal trading rate is explicitly given by
\begin{equation*}
\alpha_t^* = - \frac{ 1 }{T-t + \gamma/\eta} \Bigg\{  Y_t^*  +  \frac{1}{2 \gamma} \bigg[ - \frac{1}{\gamma_0}\sqrt{\frac{\pi}{2}}e^{-\frac{b_0^2}{2\gamma_0^2}}\left(\text{erfi}\left(\frac{b_0+\gamma_0^2
  (T-t)}{\sqrt{2}\gamma_0}\right)-\text{erfi}\left(\frac{b_0}{\sqrt{2}\gamma_0}\right)\right) +(T-t+\frac{\gamma}{\eta})  \bigg] S_t \Bigg\},
\end{equation*}
where  \texttt{erfi} is the imaginary error function, defined as:
\[
\text{erfi}(x)=\frac{2}{\sqrt{\pi}} \int_{0}^{x} e^{t^2}dt.
\]

\begin{remark}
  In the particular case where the price process is a martingale, i.e., $b_0$ $=$ $0$,
  and in the limiting case when the penalty parameter $\eta$ goes to infinity, corresponding to the final constraint $Y_T$ $=$ $0$, we see that $\alpha_t^*$ converges to
  $-Y_t^*/(T-t)$, hence it becomes independent  of the price process, and this leads to an explicit optimal inventory: ${Y_t^*=y_0 \frac{T-t}{T}}$ with constant trading rate $\alpha_t^*$ $=$
  $-y_0/T$.
  We retrieve the well-known VWAP strategy obtained in \cite{almCri}.

\end{remark}

\vspace{3mm}
\noindent We solve the problem numerically, taking $N=100$ for the time discretization, and fixing the other parameters as follows: $\gamma$=5,  $S_0$=6, $Y_0$=1, $\eta$=100 and $\sigma$=0.4.
We run two sets of forward Monte Carlo simulations for $b_0=0.1$, $T=1$ and  $b_0=-0.1$, $T=0.5$ changing the value of $\gamma_0$.
We tested the Regress-Later Monte Carlo (RLMC), the Control Randomization (CR) and the Quantization (Q) algorithms.
In particular, we wanted to compare the performance of these algorithms with $(\alpha^*_{t_n})_{n=0}^{N-1}$, where $\alpha^*$, defined above, is the optimal strategy associated to the continuous-time Portfolio Liquidation problem.
We refer to this discrete-time strategy as $\alpha^*$ (i.e., re-using the same  notation), and we use Opt, or \textit{continuous-time} optimal strategy when we want to stress the fact that this strategy is optimal for the \textit{continuous-time} control problem, and not for the discrete time one.
We also tested a benchmark strategy (Bench) which consists in liquidating the inventory at a constant rate $-y_0/T$.
The test consisted in computing the estimates $\doublehat V_{\Delta t}(t_0=0,S_0=6,Y_0=1)$ associated to the different algorithms.

\begin{figure}
  \begin{subfigure}{.48\linewidth}\centering
  {
  \includegraphics[width=1\linewidth]{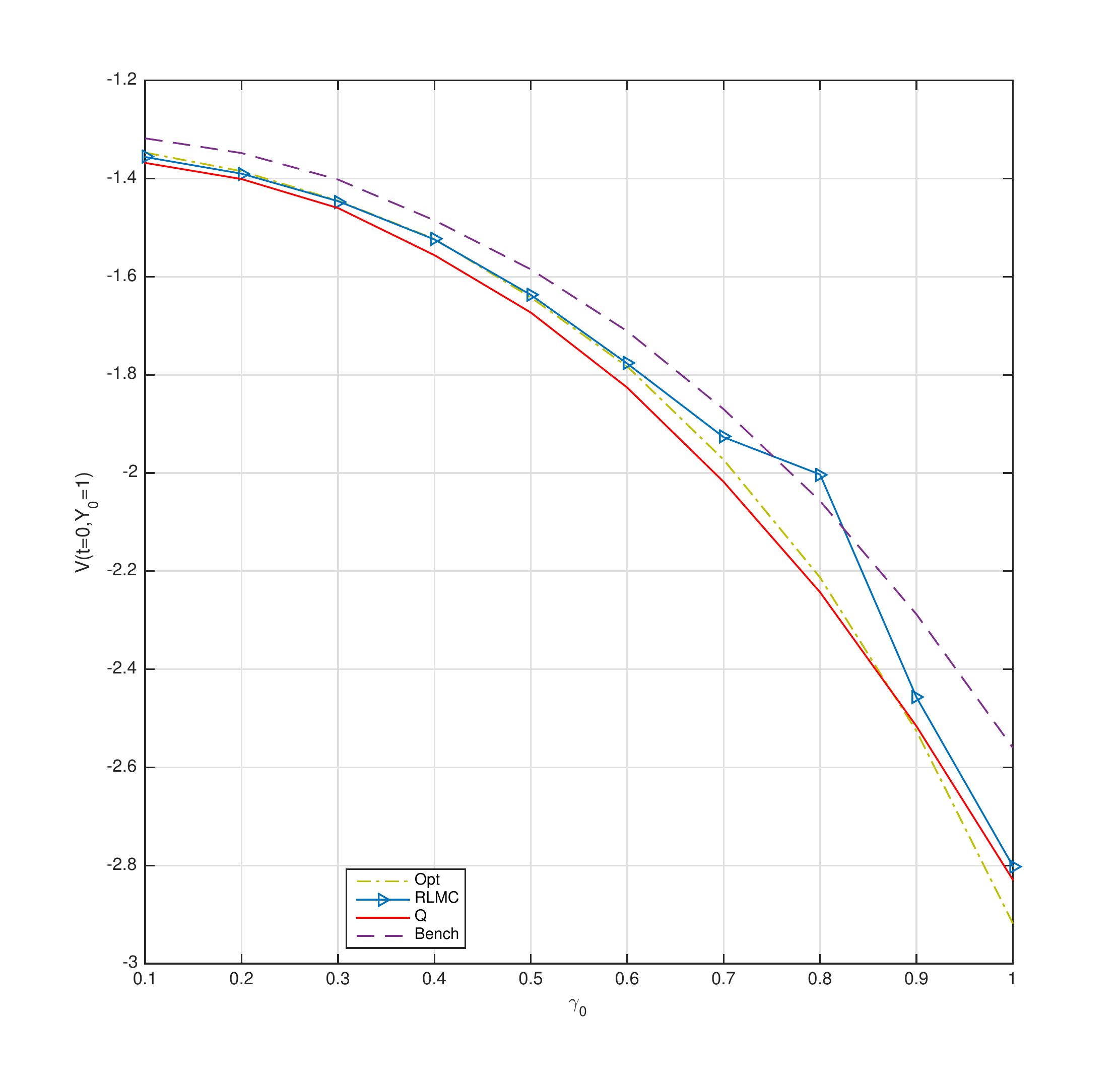}
  }
  \vspace*{-1cm}
  \caption*{$b_0=0.1$ and $T=1$.}\label{fig:PL1}
  \end{subfigure}
  \begin{subfigure}{.48\linewidth}\centering
  {
  \includegraphics[width=1\linewidth]{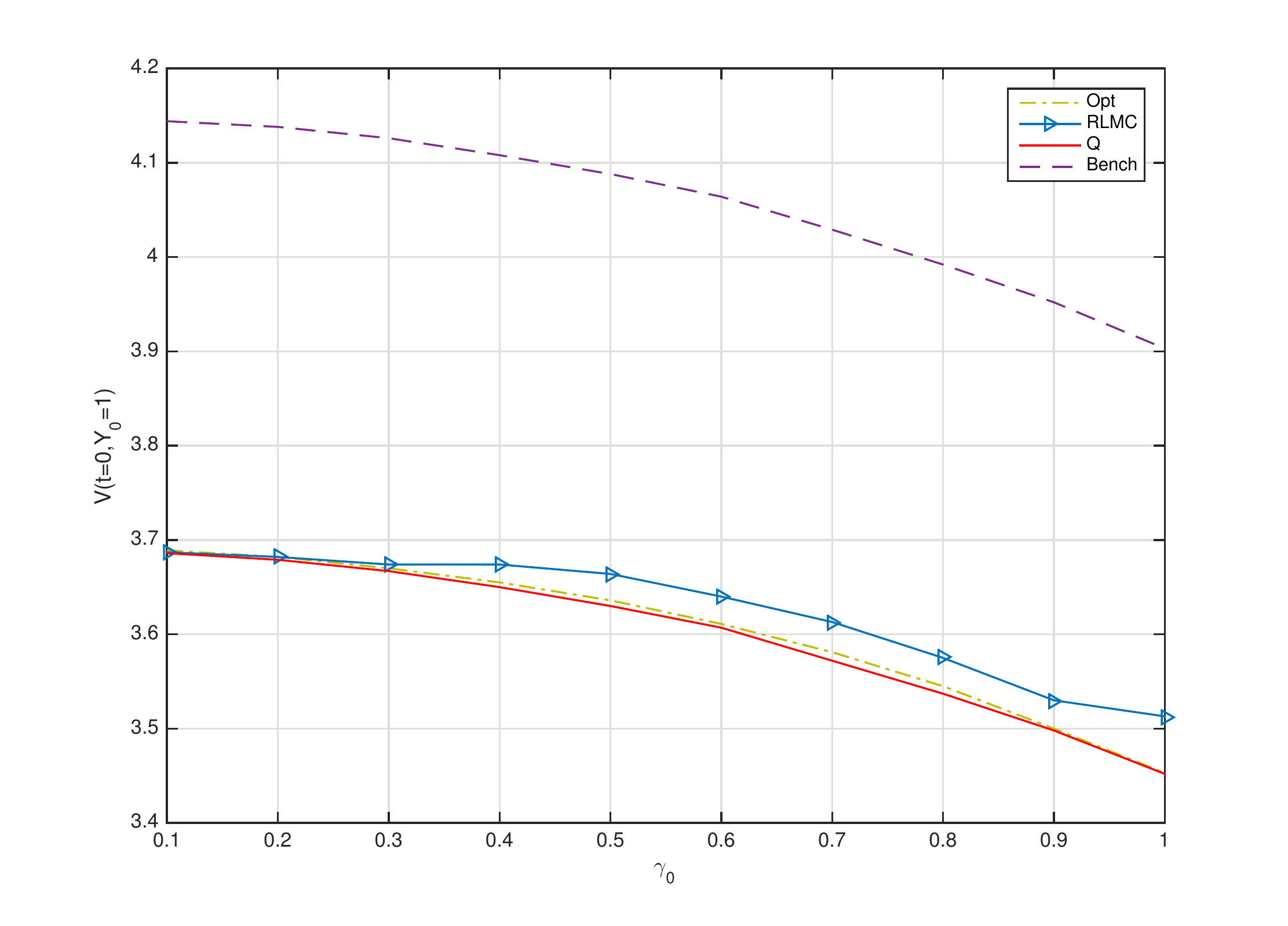}
  }
  \vspace*{-.8cm}
  \caption*{$b_0=-0.1$ and  $T=0.5$.}\label{fig:PL2}
  \end{subfigure}
  \caption{Results for the Portfolio Liquidation problem.
  Estimation of the value function at point $(s_0=6,y_0=1)$ at time 0 provided by different strategies w.r.t. $\gamma_0$. We took $\gamma$=5,  $S_0$=6, $Y_0$=1, $\eta$=100 and $\sigma$=0.4. }\label{fig:PL}
\end{figure}

\vspace{3mm}
\noindent We display the results obtained by the different algorithms in Table \ref{t:PL} and plot them in Figure \ref{fig:PL}.
One can observe in Figure \ref{fig:PL} that for $\Delta t = \frac{1}{100}$ the estimations $\doublehat V_{\Delta t}(t_0=0,S_0=6,Y_0=1)$ of the value function $V_{\Delta t}(t_0=0,S_0=6,Y_0=1)$, provided by RLMC, CR or Q-based methods, are sometimes such that
\[
\doublehat V_{\Delta t}(t_0=0,S_0=6,Y_0=1) \leq \hat J_{\Delta t}(t_0=0,S_0=6,Y_0=1,\alpha^*),
\]
where $\hat J_{\Delta t}(\cdot,\cdot,\cdot,\alpha^*)$ is a Monte Carlo estimate of $ J_{\Delta t}(\cdot,\cdot,\cdot,\alpha^*)$ applying strategy $(\alpha^*_{t_n})_{n=0}^{N-1}$ (see in Figure \ref{fig:PL} the curve Opt).
It means that RLMC, CR, or Q-based methods sometimes provide better estimations of the optimal strategy than $\alpha^{*}$ for the discrete time control problem.
However, since under suitable conditions (see, e.g., \cite{khalanpha13}), the optimal strategy for the discrete time control problem $\alpha^*_{\Delta t}$  converges toward $\alpha^*$, i.e. we have $\alpha^*_{\Delta t} \xrightarrow[\Delta t \to 0]{} \alpha^*$,  then it holds:
\[
 \hat  J_{\Delta t}(t_0=0,S_0=6,Y_0=1,\alpha^*) \xrightarrow[\Delta t \to 0]{}  V(t_0=0,S_0=6,Y_0=1).
\]

\begin{figure}
  \centering
  \includegraphics[width=.7\linewidth]{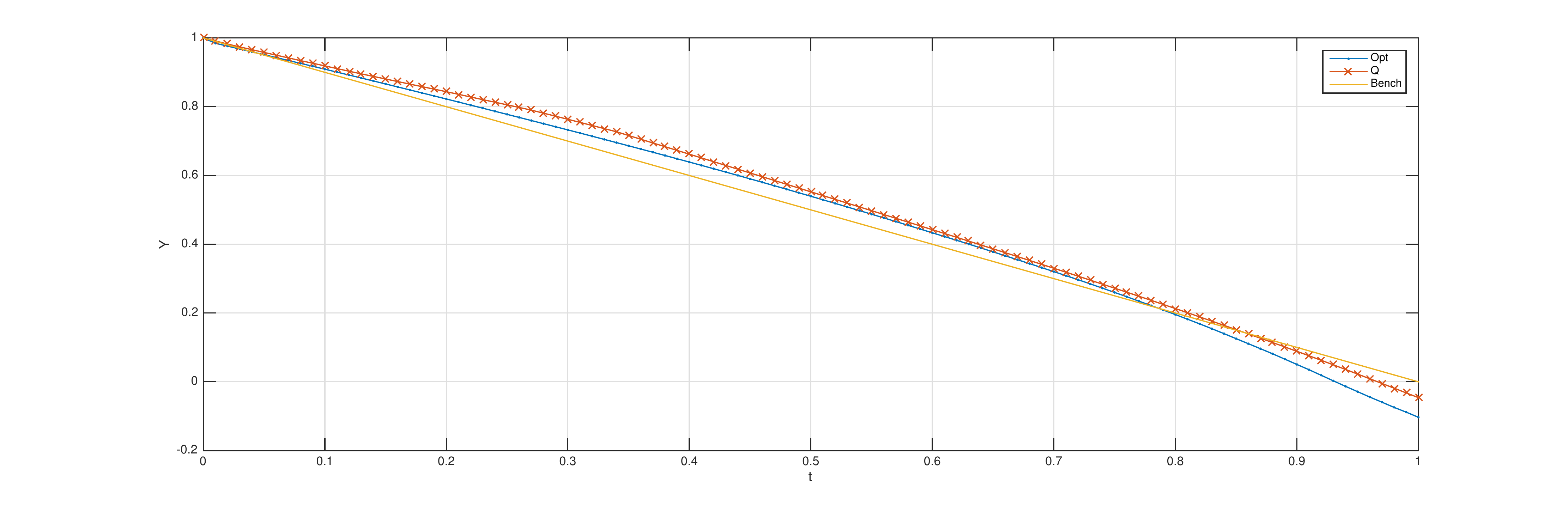}
  \caption{Simulation of $(Y_t)_{t \in [0,T]}$ using the (continuous-time) optimal strategy (Opt), the (Q) estimated one, and the Benchmark strategy (Bench) to solve the Portfolio Liquidation problem.
  We took $T=1$, $\sigma=0.4$, $\gamma_0=1$, $b_0=0.1$, $S_0=6$, $Y_0=1$, $N=100$, $\gamma=5$, $\eta=100$.
}
  \label{fig:QPL1Sample}
\end{figure}

\vspace{2mm}
Figure \ref{fig:QPL1Sample} shows a sample of the inventory $(Y_t)_{t \in [0,T]}$ when the agent follows $\alpha^*$ and the Quantization algorithm.
One can notice that given the chosen penalization parameters, it is optimal to short some stocks at terminal time.
Finally, notice that the concavity of the curves comes from the fact that the running cost does not penalize the inventory.
If so, we expect the curves of the inventory w.r.t. time to be convex, see, e.g., \cite{gasch13}.

\vspace{3mm}

\noindent {\bf Details on the RL and CR algorithms implementation} \\
The implementation of Regression Monte Carlo algorithms has required intense tuning and the use of the performance iteration technique introduced in Subsection \ref{sec:explorefirstregresslater}  in order to obtain satisfactory results.
Paramount is, in addition, the distribution chosen for the training points in Regress-Later and for the initial control in Control Randomization.
The problem of finding the best set of data to provide to the backward procedure is similar in the two Regression Monte Carlo algorithms. However little study is available in the literature; for more details on this problem in the Regress-Later setting see \cite{nadetal17} and \cite{balpal17}.
In the case of RL algorithm a training measure $\mu_n$ has been chosen in order to sufficiently explore the state space in the $Y$ dimension, in particular we considered $\mu_n = \mathcal{U}[-0.5,0.5+\frac{T-t_n}{t_n}]$. Similarly for CR we seek a distribution of the random control such that the controlled process $Y$ results in having a distribution similar to $\mu_n$.
In order to achieve such goal we follow the ``explore first, exploit later'' approach presented in Subsection \ref{sec:explorefirstregresslater} and use a perturbed version of the  empirical distribution of the control (to avoid concentration of the training points) obtained at previous iteration of the method to determine the random control at next iteration of the method.

In order to choose the basis functions, we used the fact that we expect the value function to be convex in the $Y$ dimension with minimum around the optimal inventory level and monotone in the $S$ dimension.
For RL algorithm we choose therefore the following set of basis functions: $\{s,y,y^2,sy,sy^2\}$, where we take the square function $y^2$ as a general approximator for convex functions around their minima (where we expect the measure $\mu_n$ to be concentrated).
On the other hand, CR requires that we guess what the functional form of the conditional expectation of the value function is with respect to the control process.
Considering our argument on square function approximating general convex functions we choose to add the set $\{\alpha,\alpha^2,\alpha y, \alpha s\}$ to the set of basis functions used by RL.\\
 Note that there is no need for time-consuming optimal control searching with such a choice of basis functions, as explained in Remark \ref{rq:instantaneousFinding}.\\

\noindent Finally note that we observed very high volatility in the quality of the policy estimated by control randomization. For this reason we estimated the policy $50$ times, and report in Table \ref{t:PL} the results provided by the best performing one; increasing the number of training points further affects the variability only marginally.

\vspace{3mm}
\noindent {\bf Details on the Q algorithm implementation} \\
\label{rem:QalgoPortfolioLiq}
To numerically solve this example, we used the optimal grid of the Gaussian random variable with $L=50$ points, denoted by $\Gamma^\eps_L$, to define the grid\footnote{We use the notation $t B = \{ t b, b \in B \}$, where $t \in \R$ and $B$ is a set.} $\Gamma_n^W= t_n\Gamma^\eps_L$ that discretizes $W_{t_n}$, the Brownian motion at time $t_n$, and the grid $\Gamma_n^Y= Y_0 - \frac{t_n}{T}+ t_n\Gamma_L^\eps$ that discretizes $Y_{t_n}$, the inventory at time $t_n$, $ \text{ for } n=0,...,N$.
Note that $\Gamma_n^Y$, for $n=0,...,N$, is centered at point $Y_0-\frac{t_n}{T}$ because we guessed that the optimal liquidation rate was close to $\frac{Y_0}{T}$ (see Figure \ref{fig:QPL1Sample} to check that our guess is correct). \\
We then considered the grid $\Gamma_n=\Gamma_n^W \times \Gamma_n^Y$ to discretize $Z_{t_n}=(W_{t_n},Y_{t_n})$, $n=0,...,N$.\\

\begin{figure}
   \centering
   \includegraphics[width=.5\linewidth]{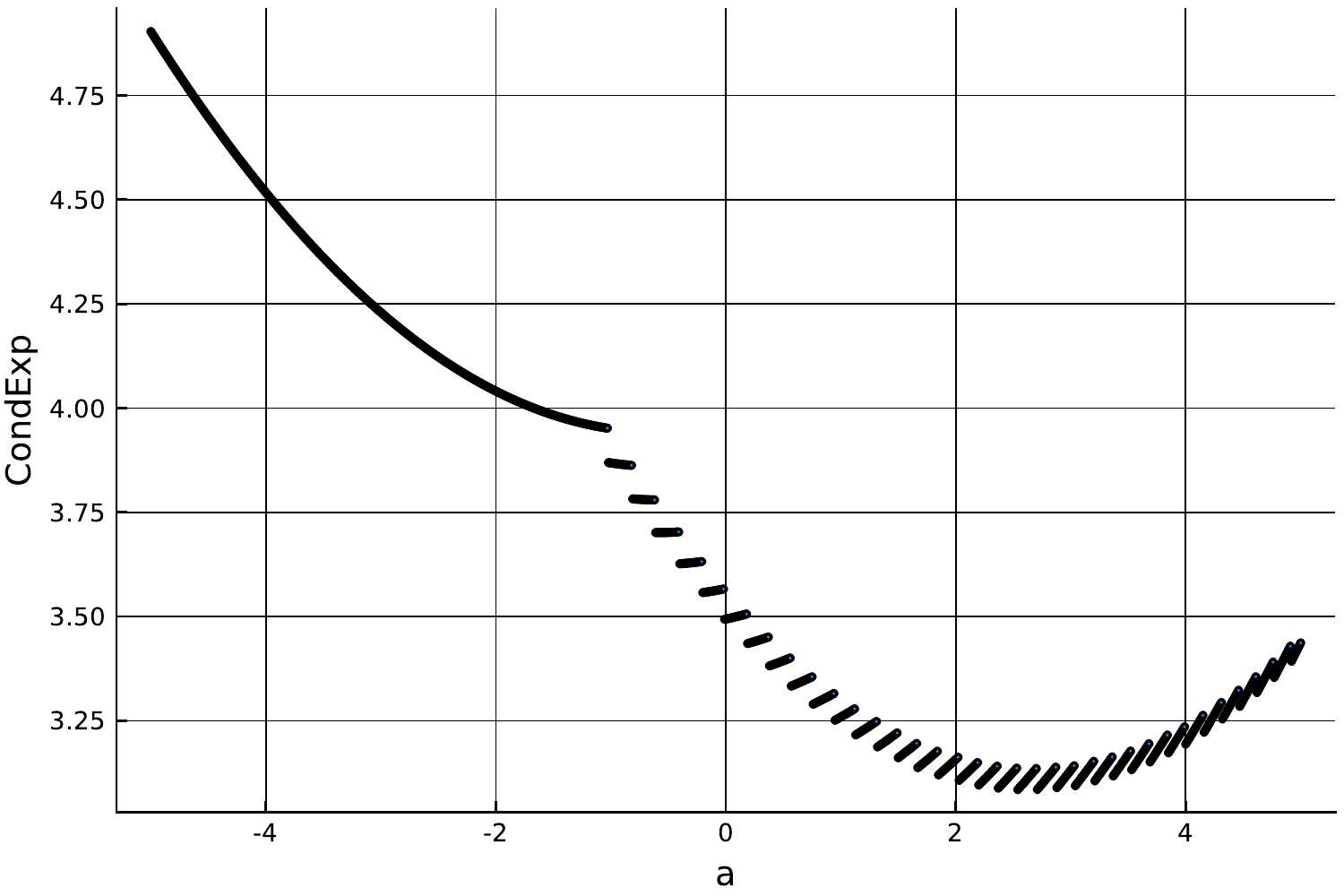}
   \caption{Plot of the quantized-based piecewise-constant approximation of the conditional expectation CondExp:\\${ a \mapsto \sum_{e \in \Gamma^\eps} \mathbb{P}(\widehat{\eps}=e)\widehat{V}_{\Delta t}^\text{Q}\Big(t_{n+1}, \text{Proj}_{\Gamma^W_{n+1}} \big(G_{\Delta t}((w,y),a,e)\big) ,\text{Proj}_{\Gamma^Y_{n+1}}  \big(G_{\Delta t}((w,y),a,e)\big)\Big)}$.\\
  We took $n=N-1$, $w=0$, and $y=-0.18$ to plot the curve.
  Observe that the approximation is discontinuous w.r.t. the control variable $a$ in such a way that it makes the search of the minimizer of this function very difficult by usual (gradient descent-based) algorithms.
  Also, observe that the minimum of the function, which is actually equal to the estimation of the value function at time $N-1$ at point $(w=0,y=-0.18)$, suffers from inaccuracy.}
   \label{fig:DiscApprox}
 \end{figure}

\vspace{2mm}
\noindent We first tried to design a quantization algorithm using the following expression for the conditional expectation approximations:
\begin{align}
\label{approx1:liquidation}
\mathbb{E}_{n,(w,y)}^a\Big[\widehat{V}_{\Delta t}^\text{Q}\Big(t_{n+1}, \text{Proj}_{\Gamma^W_{n+1}} \big( W_{t_{n+1}} \big) , \text{Proj}_{\Gamma^Y_{n+1}}\big( Y_{t_{n+1}} \big) \Big)\Big] \\
& \hspace{-5cm} \approx \sum_{\ell=1}^L p_\ell\widehat{V}_{\Delta t}^\text{Q}\Big(t_{n+1}, \text{Proj}_{\Gamma^W_{n+1}} \big(G_{\Delta t}((w,y),a,e_\ell)\big) ,\text{Proj}_{\Gamma^Y_{n+1}}  \big(G_{\Delta t}((w,y),a,e_\ell)\big)\Big), \nonumber \\
&\hspace{3cm} \text{ for } (w,y,a) \in \Gamma^W_n \times \Gamma^Y_n \times A,
\notag
\end{align}
where  the first and second components of the process $Z=(W,Y)$ are projected respectively on the grids $\Gamma^W_n$ and $\Gamma^Y_n$; and $\text{Proj}_{\Gamma^W_{n}}$  (resp. $\text{Proj}_{\Gamma^Y_{n}} $) stands for the Euclidean projection of the first (resp. second) component of $Z=(W,Y)$ on $\Gamma^W_n$ (resp. $\Gamma^Y_n$).\\
 This approximation belongs to the family of constant piecewise approximations, and is in the spirit of multidimensional component-wise-quantization methods already studied in the literature (see, e.g., \cite{PaSa}).\\
 Unfortunately, as it can be seen in Figure \ref{fig:DiscApprox}, approximation \eqref{approx1:liquidation} is discontinuous w.r.t. the control variable $a$ in such a way that the optimal control searching task suffered from instability and inaccuracy, which implied bad value function estimations at time $n=0,...,N-1$.
 We thus had to use a better conditional expectation approximation.\\

 \vspace{3mm}

\noindent We then decided to smooth the previous approximation of the conditional expectations w.r.t. the control variable by considering the following
{\begin{align*}
\mathbb{E}_{n,(w,y)}^a\Big[\widehat{V}_{\Delta t}^\text{Q}\Big(t_{n+1},  \text{Proj}_{\Gamma^W_{n+1}} \big( W_{t_{n+1}} \big) , \text{Proj}_{\Gamma^Y_{n+1}}\big( Y_{t_{n+1}} \big) \Big)\Big] &\\
& \hspace{-7cm} \approx \sum_{\ell=1}^L p_\ell \bigg[\lambda^{e_\ell,(w,y)}_{a}\widehat{V}_{\Delta t}^\text{Q}\big(t_{n+1},\text{Proj}_{\Gamma^W_{n+1}} \big[ G_{\Delta t}^w((w,y),a,e_\ell) \big] , y_+\big) \\
& \hspace{-5.5cm}+(1-\lambda^{e_\ell,(w,y)}_a)\widehat{V}_{\Delta t}^\text{Q}\Big(t_{n+1}, \text{Proj}_{\Gamma^W_{n+1}} \big[ G_{\Delta t}^w((w,y),a,e_\ell) \big], y_- \Big) \bigg],
\end{align*}
}
where, in the spirit of the semi-linear approximation presented in Subsection \ref{sec:approxexpect}, we have for all $\ell=1,...,L$:
\begin{itemize}
  \item  $G_{\Delta t}^w((w,y),a,e_\ell)$ and $G_{\Delta t}^y((w,y),a,e_\ell)$ respectively stand for the first and the second component of $G_{\Delta t}((w,y),a,e_\ell)$, i.e., $G_{\Delta t}((w,y),a,e_\ell)=\big(G_{\Delta t}^w((w,y),a,e_\ell),G_{\Delta t}^y((w,y),a,e_\ell) \big)$.
  See \eqref{def:G_deltat} for the definition of $G_{\Delta t}$.
  \item $y_-$ and $y_+$ are the two closest states in $\Gamma^Y_{n+1}$ from $ G_{\Delta t}^y((w,y),a,e_\ell) $, such that $y_- < G_{\Delta t}^y((w,y),a,e_\ell) < y_+$ if such point exists; $y_-$ and $y_+$ are equal to the closest state in $\Gamma^Y_{n+1}$ from $ G_{\Delta t}^y((w,y),a,e_\ell) $ otherwise.
  \item $\lambda^{e_\ell,(w,y)}_a = \frac{G_{\Delta t}^y((w,y),a,e_\ell)-y_-}{y_+-y_-}$ in the first case of the definition of $y_-$ and $y_+$ above; $\lambda^{e_\ell,(w,y)}_a = 1$ otherwise.
\end{itemize}
This approximation is a slight generalization (to dimension $d$=2) of the semi-linear approximation developed in \eqref{approxEspCond2}.
Its main interest lies in the continuity of the approximation w.r.t. the control variable $a$, which provides stability and accuracy to the usual (gradient descent-based) algorithms for the optimal controls searching, as can be seen on the numerical results (see, e.g., Table \ref{t:PL}).

 \begin{table}[]
   \centering
    \caption{Portfolio Liquidation results. Estimations of the value functions at point $(s_0=6,y_0=1)$ and time 0 provided by different algorithms.}
    \label{t:PL}
   \begin{tabular}{l|lllll|lllll}
     & \multicolumn{5}{c|}{$b_0=0.1$, $T=1$} & \multicolumn{5}{c}{$b_0=-0.1$, $T=1/2$} \\ \cline{2-11}
     $\gamma_0$ & Opt     & RLMC       & CR     & Q & Bench   & Opt      & RLMC       & CR      & Q & Bench \\ \hline
      0.1   &   -1.347   &  -1.356    &   -1.278   &  -1.368   & -1.318   & 3.689 &    3.687& 3.995   & 3.686      & 4.144 \\
     0.2   &   -1.385    &  -1.390    &   -1.283   &    -1.401 & -1.348   & 3.682       &    3.682& 3.847   &   3.679    &  4.138 \\
     0.3   &  -1.445     &   -1.446   &   -1.314   &   -1.460  &  -1.402 &3.670       &    3.674& 4.034   &     3.667    & 4.126 \\
     0.4   &   -1.523    &   -1.524   &   -1.323   &     -1.556 & -1.485  & 3.655        &    3.674 &4.128    &     3.650    &  4.108 \\
     0.5   &   -1.642    &  -1.637    &  -1.348    &   -1.673   & -1.585  & 3.636        &    3.664 &4.243    &     3.630    & 4.088 \\
     0.6   &    -1.783   &  -1.777    &   -1.425   &   -1.826  &  -1.711  & 3.611        &    3.640 &4.386    &     3.607    &  4.064 \\
     0.7   &    -1.973   &  -1.927    &   -1.513   &  -2.018   &   -1.870 &3.581        &    3.613 &4.783    &     3.572    & 4.029 \\
     0.8   &   -2.213    &  -2.003    &   -1.637   &  -2.243   &  -2.057  & 3.545        &    3.575 &5.142    &     3.537    & 3.992 \\
     0.9   &   -2.526    &   -2.457   &   -1.819   &  -2.516    & -2.288  &3.500        &    3.530 &5.345   &     3.498    &  3.952 \\
     1      &    -2.918  &   -2.801   &   -1.806  &   -2.829  &  -2.560   & 3.453        &    3.513 &6.765    &    3.452     & 3.903
   \end{tabular}
 \end{table}

\vspace{10mm}

\noindent {\bf 2. Portfolio Selection.}
Consider the Portfolio Selection problem with one risky asset.
We choose a CARA utility function $U(x)$ $=$ $-\exp(-px)$, with $p$ $>$ $0$.
It has been shown in \cite[Corollary 1]{guepu16} that the optimal portfolio strategy is explicitly given by
\beqs
\alpha_t^* &=& \frac{\sigma^2 + \gamma_0^2 t }{\sigma^2 + \gamma_0^2 T} \frac{\hat\beta_t}{p\sigma^2}
\enqs
where
\beqs
\hat\beta_t &=& \E^0[ \beta | \Fc_t^S] \; = \;  \frac{\sigma^2}{\sigma^2 + \gamma_0^2 t} b_0 + \frac{\gamma_0^2}{\sigma^2 + \gamma_0^2 t} \Big(  \ln\frac{S_t}{S_0}
+ \frac{1}{2}\sigma^2 t  \Big),
\enqs
is the posterior mean of the drift  (Bayesian learning on the drift), and the optimal performance by
\beqs
 J_2(\alpha^*) &=& -  \exp\Big[  - p \Big(
 x_0 + \frac{1}{2p} \big(   \ln \big(\frac{\sigma^2 + \gamma_0^2 T }{\sigma^2}\big) -   \frac{\gamma_0^2 T}{\sigma^2 + \gamma_0^2 T}  \big) +
 \frac{b_0^2}{2p\sigma^2} \frac{\sigma^2 T}{\sigma^2 + \gamma_0^2 T}
 \Big)
  \Big].
\enqs

\vspace{3mm}
 The Portfolio Selection problem, even though in many aspects similar to the Portfolio Liquidation problem, is interesting in its own right because the control acts only on the variance of the controlled wealth process.
 We tested the Regress-Later Monte Carlo  (RLMC), the Control Randomization  (CR) and the Quantization (Q) algorithm on the Portfolio Selection problem.
 Similarly to what has been done for Portfolio Liquidation problem, we discretized time choosing $N=100$ and solved the discrete time problem associated.
 We considered two set of experiments, $b_0=0.1$, $T=1$ and $b_0=-0.1$, $T=0.5$, for different values of $\gamma_0 \in [0,1]$, $p=1$, $\sigma = 0.4$.
 Given all these different parameters, we compared the performance of these algorithms with the one of the optimal strategy for the continuous-time problem $\alpha^*$ (Opt).
 The general test consists in computing a forward Monte Carlo with 500000 samples, following optimal strategy estimated using different strategies, to provide estimates of $V(t_0=0,X_0=0,W_0=0)$ the value function at time 0.

 \vspace{3mm}
 \noindent  We present the results of our numerical experiments in Table \ref{t:PS}.
 One can see that the Quantization algorithm performs similarly to the theoretical optimal strategy (Opt) for the continuous time problem, which can be interpreted as stability and accuracy of the Q algorithm, and also shows that the time discretization error is almost zero here.

\begin{figure}
   \centering
   \includegraphics[width=1\linewidth]{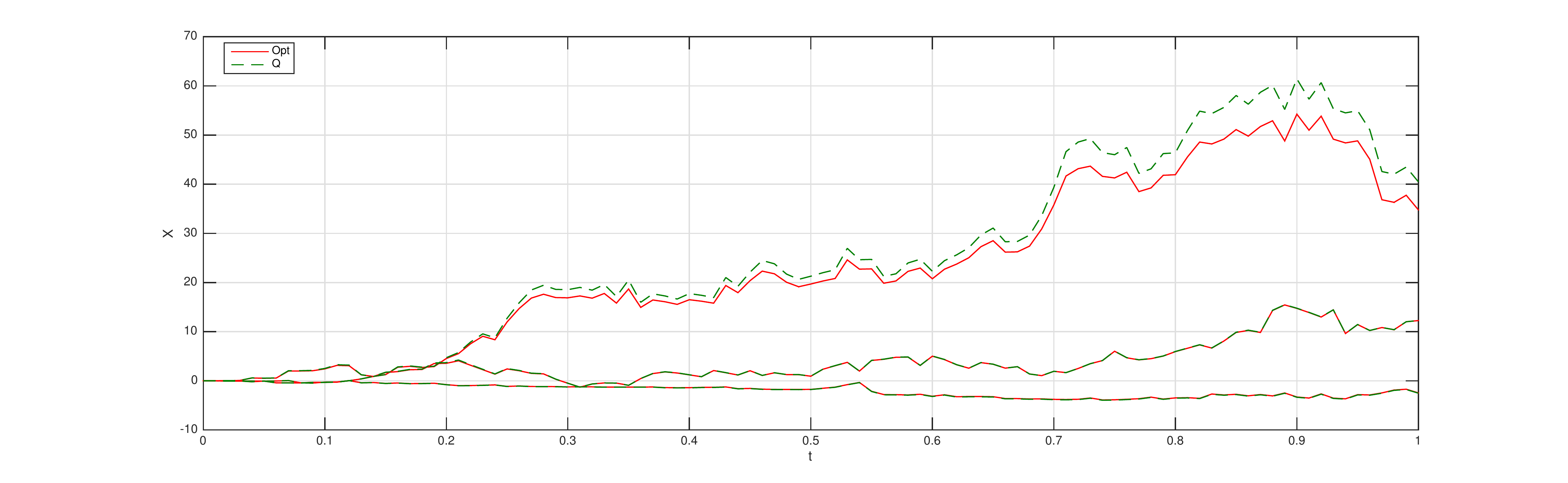}
   \caption{3 simulations of the agent's wealth $(X_t)_{t\in [0,T]}$ when the latter follows the continuous-time optimal strategy (Opt) and the (Q) estimated optimal strategy to solve the Portfolio Selection problem.
  We took $\sigma$=0.4, $T$=1, $P$=0.1, $\gamma_0$=5, $b_0$=0,1.
  One can see that the two strategies are the same when the drift is low; but Q performs slightly better than Opt when the drift is high, which is a time-discretization effect.}
   \label{fig:PSPaths}
 \end{figure}

 \vspace{3mm}
\noindent We also present in Figure \ref{fig:PSPaths} a sample of the wealth of the agent following the optimal strategy and the (Q) estimated one.
One can see that the strategies slightly differ when the drift is high, and remain the same when the drift is low.
The small difference can be explained by the fact that the optimal strategy (Opt) is not optimal for the discrete time version of the problem.

\vspace{3mm}
\noindent {\bf Details on the Q algorithm implementation} \\
We designed the same Quantization algorithm as the one built to solve the Portfolio Liquidation problem.
We nevertheless had to take a larger number of points in the grids to minimize the back-propagation of errors from the borders of the girds; and had to use the ``explore first, exploit later'' idea (see Subsection \ref{sec:explorefirstregresslater}) to improve the results.

\vspace{3mm}
\noindent {\bf Details on the RL and CR algorithms implementation} \\
When implementing Regression Monte Carlo algorithms, and choosing basis functions, the control on variance implies that low order polynomial can not be used alone, as they can easily cause the control to be bang-bang between the boundaries of its domain.
Similarly, piecewise approximations are not very effective, as the dependence on the control is very weak, requiring a high number of local supports and making the computational complexity overwhelming.
We tested both value and performance iteration and tried to employ different kinds of basis functions and training points. Unfortunately, both Regress-Later and Control Randomization do not cope well with controlling the dynamics of a process through the variance only.
A tailor-made implementation of Regression Monte Carlo to deal with this kind of problems is outside the scope of this paper and further investigation will follow in future work. For now, we chose not to provide results based on RL and CR methods.


 \begin{table}[]
	\centering
	\caption{Portfolio Selection results. Estimations of the value function at point $(x_0=0,S_0=6)$ time 0 using the \textit{continuous-time} optimal strategy (Opt) and (Q) estimated optimal strategy.}
	\label{t:PS}
	\begin{tabular}{l|ll|ll}
		& \multicolumn{2}{c|}{$b_0=0.1$, $T=1$} & \multicolumn{2}{c}{$b_0=-0.1$, $T=0.5$} \\ \cline{2-5}
		$\gamma_0$ & Opt        & Q    & Opt       & Q      \\ \hline
		0.1        &   -0.985  &    -0.985  &    -0.992   &    -0.992    \\
		0.2        &  -0.982    &     -0.982 &  -0.991       &    -0.991   \\
		0.3        &  -0.973    &  -0.973    &  -0.988      &    -0.988   \\
		0.4        &  -0.954       &    -0.953  &      -0.981 &       -0.981 \\
		0.5        &    -0.927    &    -0.927  &   -0.969    &     -0.969   \\
		0.6        &   -0.896    &    -0.896    &     -0.952       &   -0.952   \\
		0.7        &  -0.863        &    -0.863  &       -0.932    &    -0.932    \\
		0.8        &    -0.830       &  -0.830    & -0.910          &     -0.910   \\
		0.9        &    -0.797       &    -0.797  &    -0.886       &    -0.886    \\
		1        &   -0.767       &   -0.766  &   -0.863      &   -0.863     \\
	\end{tabular}
\end{table}

\subsection{A model of interbank systemic risk with partial observation}

\subsubsection{The model}

We consider the following model of systemic risk inspired by the model in \cite{caretal14}. The log-monetary reserves of $N$ banks lending to and borrowing from each other are governed by the system
\beqs
dX_t^i &=& \frac{\kappa}{N} \sum_{j=1}^N (X_t^j-X_t^i) dt + \sigma X_t^i ( \sqrt{1-\rho^2} dW_t^i + \rho  dW_t^0)  , \;\;\; i=1,\ldots,N
\enqs
where $W^i$, $i$ $=$ $1,\ldots,N$, are independent Brownian motions, representing the idiosyncratic risk of each bank, $W^0$ is a common noise independent of $W^i$, $\sigma>0$ is given real parameter, $\rho$ $\in$ $[-1,1]$, and where $X^i_0$, $i=1,...,N$ are i.i.d..
The mean-reversion coefficient $\kappa$ $>$ $0$ models the strength of interaction between the banks where bank $i$ can lend to and borrow from bank $j$ with an amount proportional to the difference between their reserves.
In the asymptotic regime when $N$ $\rightarrow$ $\infty$, the theory of propagation of chaos implies that the
reserve state $X^i$  of individual banks become independent and identically distributed conditionally on the common noise $W^0$, with a state governed by
\beqs
dX_t &=& \kappa (\E[X_t|W^0] - X_t) dt + \sigma X_t ( \sqrt{1-\rho^2} dB_t + \rho  dW_t^0)
\enqs
for some Brownian motion $B$ independent of $W^0$.

Let us now consider a central bank, viewed as a social planner, who only observes the common noise and not the reserves of each bank, and can influence the strength of the interaction between the individual banks, through an $\F^0$-adapted control process $\alpha_t$.
The reserve of the representative bank in the asymptotic regime is then driven by
\beqs
dX_t &=& (\kappa + \alpha_t) (\E[X_t|W^0] - X_t) dt + \sigma X_t ( \sqrt{1-\rho^2} dB_t + \rho  dW_t^0),\quad X_0 \sim X^1_0,
\enqs
and we consider that the objective of the central bank is to minimize
\beqs
J(\alpha) &=& \E \Big[ \int_0^T \left( \frac{1}{2}\alpha_t^2 + \frac{\eta}{2} (X_t - \E[X_t|W^0])^2 \right) dt + \frac{c}{2}(X_T- \E[X_T|W^0])^2 \Big],
\enqs
where $\eta$ $>$ $0$ and $c$ $>$ $0$ penalize the departure of the reserve from the average.
This is a
MKV control problem under partial observation, but notice that it does not belong to the class of linear quadratic (LQ) MKV problems due to the control $\alpha$ which appears in a  multiplicative form with the state.
However, it fits into our class of polynomial MKV problem, and can be embedded into standard control problem  as follows:  We set $\bar X_t$ $=$ $\E[X_t|W^0]$ and $Y_t$ $=$ $\E[(X_t-\bar X_t)^2|W^0]$.  The cost functional is then written as
\beqs
\label{eq:HJBSysRisk}
J(\alpha) &=& \E \Big[ \int_0^T \left( \frac{1}{2}\alpha_t^2 +  \frac{\eta}{2} Y_t \right) dt +  \frac{c}{2} Y_T \Big]
\enqs
where the dynamics of $\bar X$ and $Y$ are governed by
\beqs
\label{eq:euler_sys}
d\bar X_t &=& \sigma \rho  \bar X_t dW_t^0, \quad \bar{X}_0= x_0 = \E[X_0] \\
dY_t &=& \big[ \big(\sigma^2 - 2(\kappa +\alpha_t) \big) Y_t + \sigma^2(1-\rho^2) \bar X_t^2 \big] dt + 2 \rho\sigma Y_t dW_t^0, \quad Y_0=\mathrm{Var}(X_0).
\enqs
We have then reduced the problem to a $(\P,\F^0)$-control problem in dimension two with state variables $(\bar X,Y)$, which is neither LQ, but can be solved numerically.

\vspace{3mm}
\subsubsection{Numerical results}

\noindent For this problem, in the absence of analytical solution, we decided to compare the estimations of the value function at time 0 provided by our algorithms with a numerical approximation based on finite difference scheme provided by Mathematica, of the solution to the $2$-dimensional HJB equation associated to the systemic risk problem:
\begin{equation}
\label{HJB:systemicRisk}
  \begin{cases}
  \partial_t V + \frac{\eta}{2}y +\Big( \big( \sigma^2-2\kappa  \big) y + \sigma^2(1-\rho^2)x^2\Big)\partial_y V+   \Sup_{a \in A} \left[ \frac{1}{2} a^2 -2ay \partial_y V  \right] \\
  \hspace{3.3cm}+ \frac{\sigma^2\rho^2x^2}{2} \partial_{xx}^2 V + 2 \sigma^2\rho^2xy \partial_{xy}V + 2 \sigma^2\rho^2 y^2 \partial_{yy}^2V =0, \quad \text{ for } (t,x,y) \in [0,T) \times \R \times \R_+, \\
  V(T,x,y)= \frac{c}{2}y, \quad \forall (x,y) \in \R \times \R_+.
  \end{cases}
\end{equation}
We refer to the solution of this partial differential equation (obtained using Mathematica using finite differences as explained below) as the Benchmark (or simply Bench) in the sequel.\\

\noindent We computed $\doublehat{V}_{\Delta t}(t_0=0,x_0=10,y_0=0)$ using RL, CR and Q methods by considering a sample of size 500 000, and using the following parameters $T=1$,    $\sigma= 0.1$,    $\kappa= 0.5$ and $X_0 = 10$.
We recall that $\doublehat{V}_{\Delta t}(t_0=0,x_0,y_0)$  is an estimation of $V(0,x_0=10,y_0=0)$, the value function at $(x_0,y_0)$ and time $0$ (see its definition on the last step of each pseudo-code presented in Subsection \ref{sec:pseudocodes}).

\begin{figure}
  \begin{subfigure}{.5\linewidth}\centering
  {
  \includegraphics[width=1\linewidth]{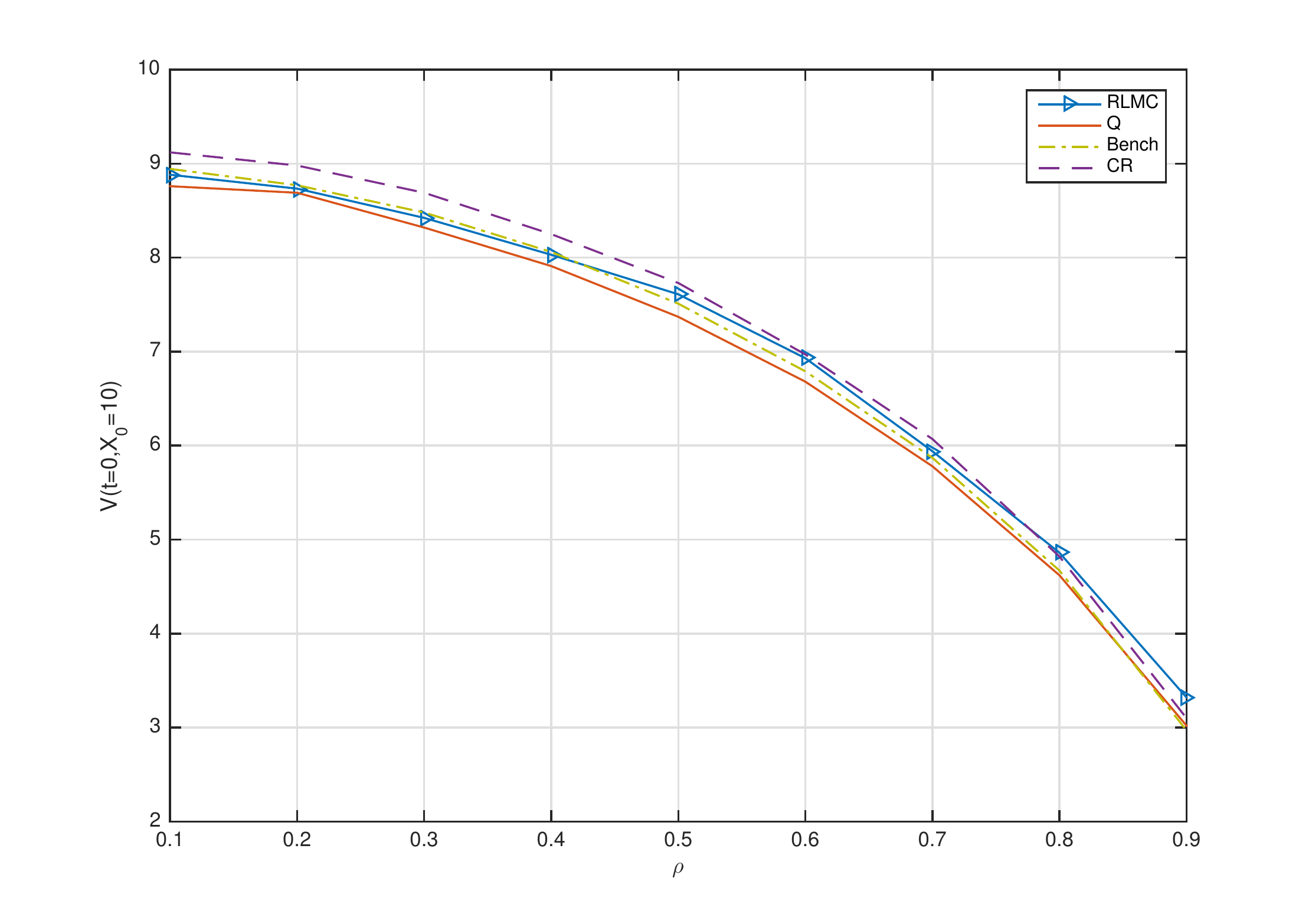}
  }
\vspace*{-0.8cm}
  \caption*{$c=100$ and $\eta=10$.}\label{fig:rho}
  \end{subfigure}
  \begin{subfigure}{.5\linewidth}\centering
  {
  \includegraphics[width=1\linewidth]{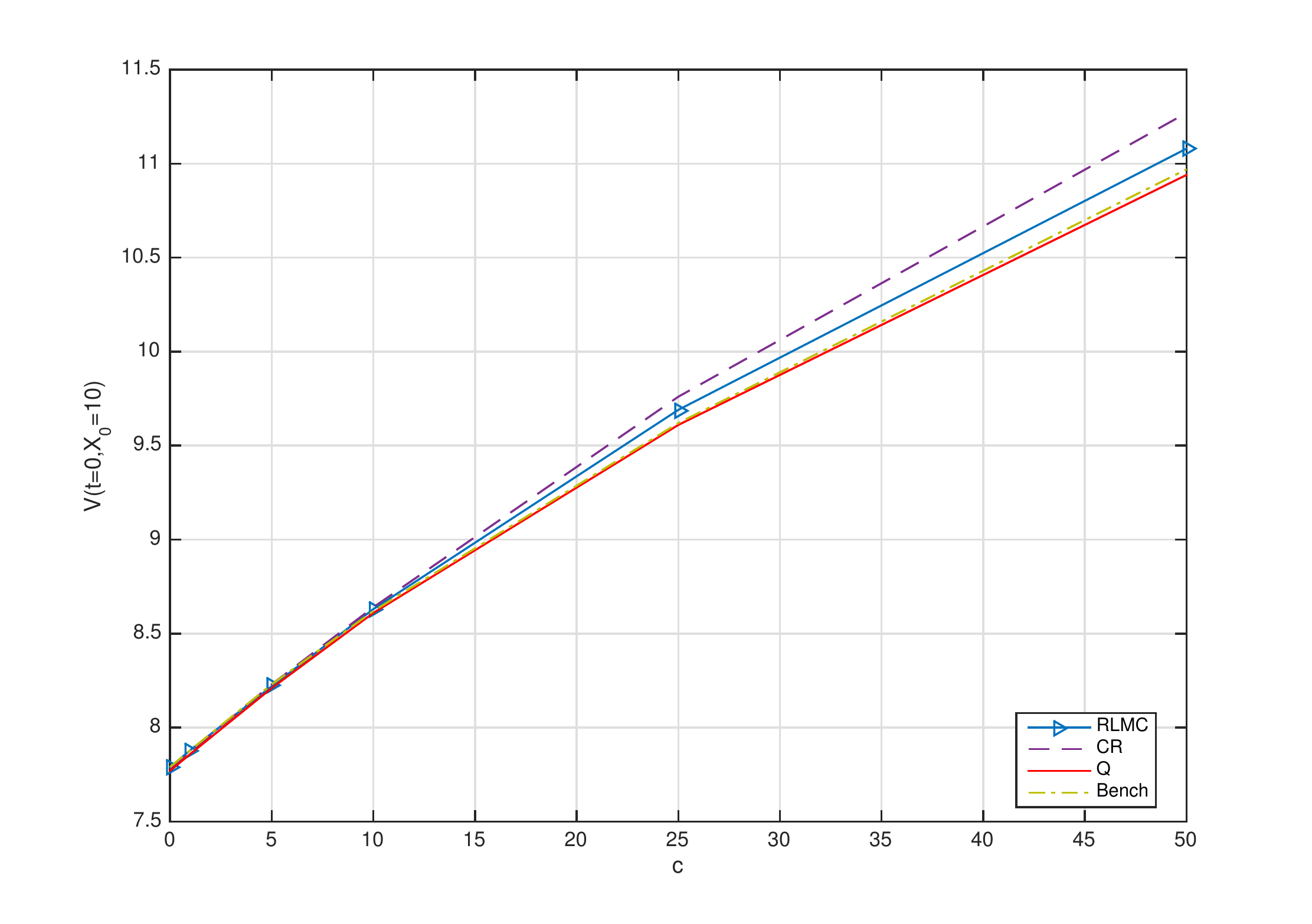}
  }
\vspace*{-0.8cm}
  \caption*{$\rho=0.5$ and  $\eta=100$.}\label{fig:c}
  \end{subfigure}
  \caption{Results for the systemic risk problem.
  Estimations of the value function at time 0 using different algorithms w.r.t. $\rho$ and $c$.
  We took $T$=1, $N$=100, $\sigma$=0.1, $\kappa$=0.5, $x_0$=10. }\label{fig:SR}
\end{figure}

\vspace{3mm}
\noindent In Table \ref{t:SR} we display the numerical results of experiments run for two situations: we took $\eta=10$, $c=100$ and  $\eta=100$, $\rho=0.5$ and vary the value of $\rho$ in the first case, and vary the value of $c$ in the second one.
Plots of the two tables are also available in Figure \ref{fig:SR}.
One can observe that the algorithms performs well.
Mainly, Bench and Q provide slightly better results than the Regression Monte Carlo-based algorithms (the curves of Bench and Q are below those of the other two).

\begin{figure}
  \centering
  \includegraphics[width=1\linewidth]{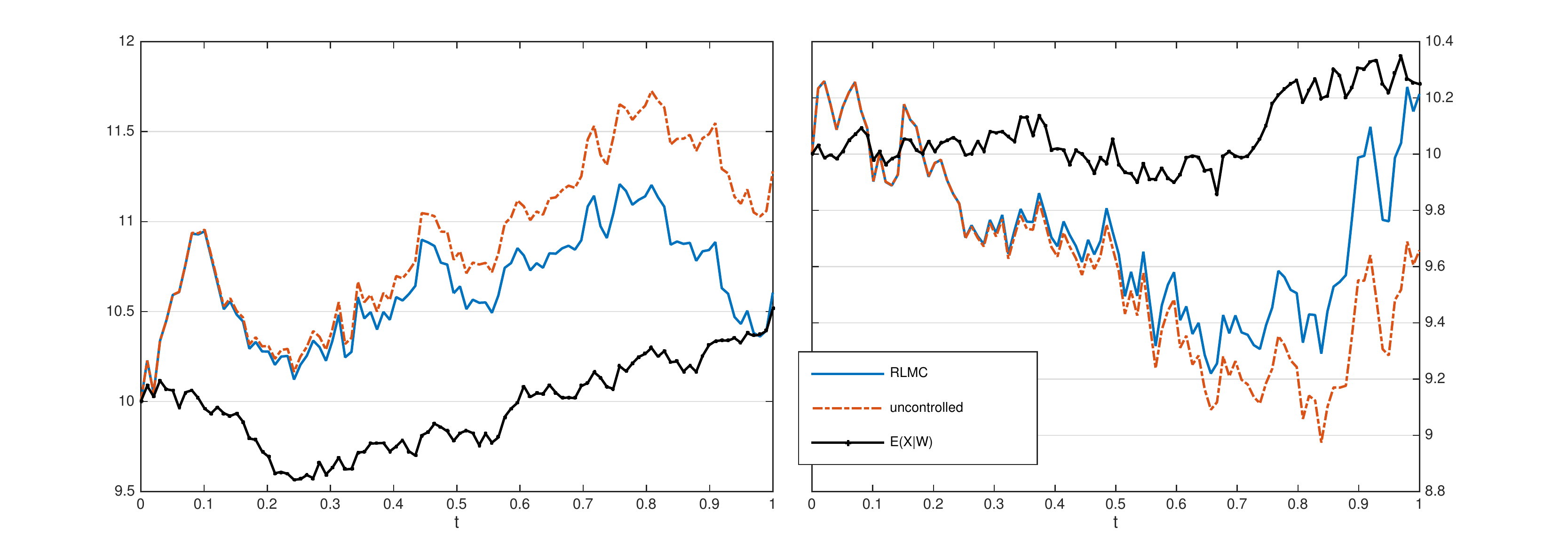}
  \caption{Two realizations of $(X_t)_{t \in [0,T]}$ controlled by RLMC (curve ``RLMC''), $(X_t)_{t \in [0,T]}$ naively controlled taken $\alpha=0$ (curve ``uncontrolled''), and $\bar{X}$ (curve ``$E(X \vert W)$'').  The optimal control for the systemic risk problem (computed by RLMC) is to do nothing at first, and catch $\bar{X}$ when the terminal time is getting close.
  }
  \label{fig:SR_pathwise}
\end{figure}

\vspace{3mm}
 Figure \ref{fig:SR_pathwise} shows two examples of paths $(X_t)_{t \in [0,T]}$ controlled by RLMC (curve ``RLMC''), $(X_t)_{t \in [0,T]}$ naively controlled by $\alpha=0$ (curve ``uncontrolled''), and the conditional expectation of $X$  $(\bar{X}_t)_{t \in [0,T]}$ (curve ``$E(X \vert W)$'').  One can see in these two examples that the (RLMC estimated) optimal control is as follows:
 \begin{itemize}
  \item  do nothing when the terminal time is far, i.e., take $\alpha=0$, not to pay any running cost.
  \item  catch $\bar{X}$ when the terminal time is getting close, to minimize the terminal cost.
 \end{itemize}

\begin{figure}
  \centering
  \includegraphics[width=1\linewidth]{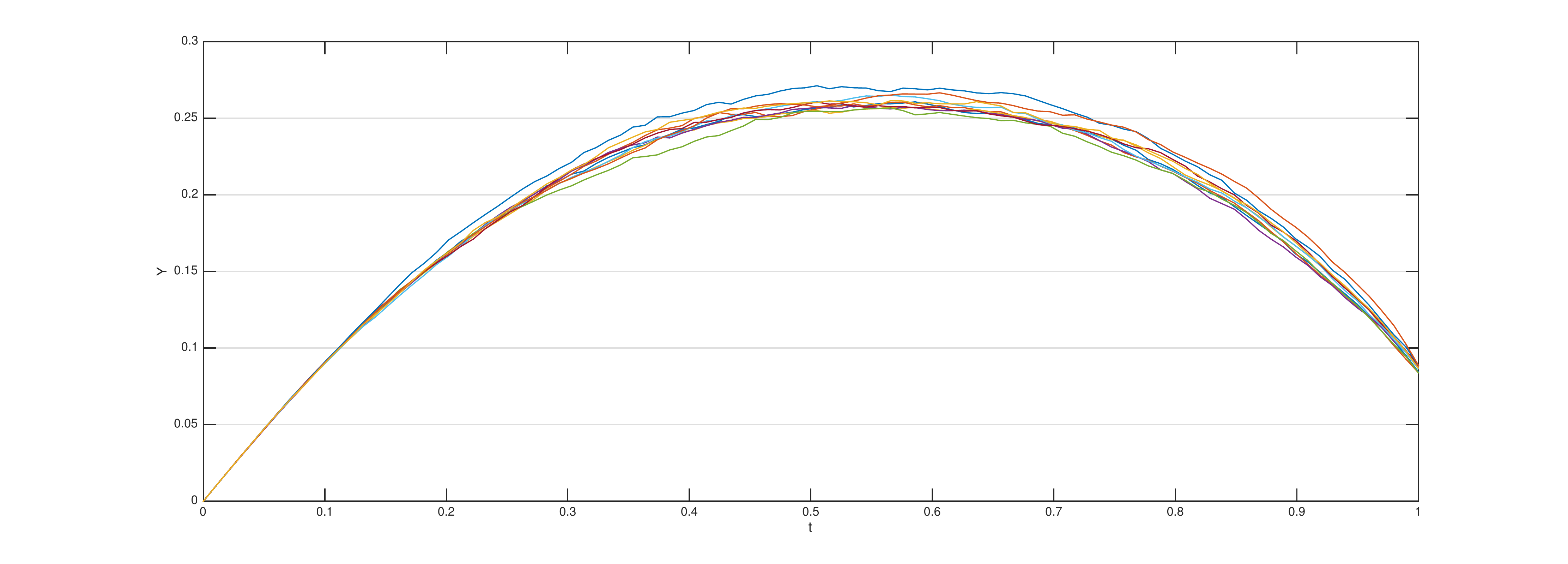}
  \caption{Sample of $(Y_t)_{t \in [0,T]}$ controlled by Q.
  The (Q) estimated optimal control for the systemic risk problem is to initially let $Y$ become large, and then reduce its value when the approaching the terminal time.
  }
  \label{fig:SR_Y}
\end{figure}

\noindent We finally present a sample of paths $(Y_t)_{t \in [0,T]}$ controlled by the decisions given by Q in Figure \ref{fig:SR_Y}.
One can see that the  (Q estimated) optimal strategy minimizes the running cost first by letting $Y$ grow; and deals with the terminal cost later by making $Y$ small when the terminal time is approaching.

\vspace{3mm}
\noindent {\bf Details on the RL and CR algorithms implementation}

For the implementation of the RL algorithm we decided to use polynomial basis functions up to degree 2.
This choice allows us to compute the optimal control analytically as a function of the regression coefficients (see Remark \ref{rq:instantaneousFinding}).
Compared to other optimization techniques, explicit expression allows for much faster and error-free computations (see Remark \ref{rq:instantaneousFinding}).
For CR, we used basis functions up to degree 3 in all dimensions to obtain more stable results.

Regarding the choice of the training measure in RL, we employed marginal normal distributions on each dimension.
As we know that the inventory dimension $Y$ represents the conditional variance of the original process $X$, we centered the training distribution $\mu_n$ at zero but considered only training points $Y_n^m\ge0$.
In CR, on the other hand, we need to carefully choose the distribution of the random control so that the process $Y$ does not become negative.
Notice in fact that the Euler approximation, contrary to the original SDE describing $Y$, does not remain positive and we would therefore need to carefully choose a control to avoid driving $Y$ negative.
In order to achieve such goal, without having to worry too much about the control, we modified the Euler approximation of \eqref{eq:euler_sys} to feature a reflexive boundary at zero. Such features allow to train the estimated control policy to not overshoot when trying to drive the process $Y$ to zero, without having $Y$ to become negative.

\vspace{3mm}
\noindent {\bf Details on the Q algorithm implementation}

As stated above, it is straightforward that $Y>0$ on $(0,T]$.
However, the Euler scheme used to approximate the dynamics of $Y$ does not prevent the associated process $(Y_{t_i})_{0<i\leq N}$ to be non-positive.
When implementing the Q algorithm for the systemic risk problem, we forced $\big(\text{Proj}_{\Gamma^Y_i}\big( Y_{t_i}\big) \big)_{0<i\leq N}$ to remain positive by simply choosing positive points for the grids $\Gamma^{Y}_i$ that quantize the states of $Y_{t_i}$, at time $t_i$ for $i=0,..., N$.

Also, given the expression of the instantaneous and terminal reward, one can expect $Y$ to stay close to 0, but we do not have any idea of how small $Y$ should stay for the strategy to be optimal (cf. Figure \ref{fig:SR_Y} to see a posteriori where $Y$ lies).
To deal with this situation, we decided to adopt the ``explore first, exploit later'' procedure.
First, we chose some random grids with a lot of points near 0 and computed the optimal strategy on these grids.
Then, we ran forward Monte Carlo simulations and generated an empirical distribution of the quantized $Y$.
Second, we build new grids of Quantization for $Y$ by generating new points according to the empirical distribution that we got from in the previous step.
Finally, we computed new (hopefully better) estimations of the optimal strategy by running the Q algorithm using the new grids.
The Q strategy performed better after applying this step, but not significantly since our first naive guess for the grids (i.e., before bootstrapping) was already good enough.

 \vspace{3mm}
\noindent {\bf Details on the implementation of the deterministic algorithm for the resolution of the HJB}

We use the \texttt{NDSolve} function in Mathematica based on finite difference method to solve \eqref{HJB:systemicRisk}.
Note that usually terminal and boundary conditions are required to get numerical results.
The final condition: \(V(T,x,y)=\frac{c}{2} y\) is already given by \eqref{HJB:systemicRisk}.
However, the boundary conditions on \(V(t,0,y)\) and \(V(t,x,0)\) are missing, except the trivial condition consisting of \(V(t,0,0) = 0\).
We then provided the HJB without boundary conditions to the Mathematica function \texttt{NDSolve}, and let the latter add artificial boundary conditions by itself to output results.

\begin{table}
  \begin{subtable}{.48\linewidth}\centering
  {\begin{tabular}{l|llll}
  \multicolumn{1}{c|}{$\rho$} & \multicolumn{1}{c}{RLMC} & \multicolumn{1}{c}{CR} & \multicolumn{1}{c}{Q} & \multicolumn{1}{c}{Bench}  \\ \hline
  0.1        &  8.88&     9.12     &  8.76 &8.94\\
  0.2       & 8.73 & 8.98 &   8.69  & 8.77\\
  0.3        & 8.42 & 8.69 &  8.32  & 8.48\\
  0.4      & 8.02 & 8.25 &  7.91   & 8.06\\
  0.5         & 7.61 & 7.73 &   7.37& 7.51\\
  0.6       & 6.93 & 6.97 &   6.68  & 6.79\\
  0.7        &  5.94 & 6.07 &  5.78 & 5.87\\
  0.8        & 4.86 & 4.82 &   4.62 & 4.67\\
  0.9        & 3.32 & 3.10  &   3.02 & 2.97
  \end{tabular}}
  \caption*{$c=100$ and $\eta=10$.}\label{tab:c}
  \end{subtable}
  \begin{subtable}{.48\linewidth}\centering
  {  \begin{tabular}{l|llll}
  \multicolumn{1}{c|}{$c$} & \multicolumn{1}{c}{RLMC} & \multicolumn{1}{c}{CR} & \multicolumn{1}{c}{Q} & \multicolumn{1}{c}{Bench} \\ \hline
  0                        & 7.79                    &        7.78                & 7.77                & 7.79                    \\
  1                        & 7.88                    &        7.87                & 7.86                & 7.88                   \\
  5                        & 8.22                    &         8.23               & 8.21                 & 8.23                     \\
  10                       & 8.63                    &         8.64              & 8.61                  & 8.62                     \\
  25                       & 9.69                    &          9.76              & 9.61                 & 9.62                   \\
  50                       & 11.08                   &           11.27             & 10.94                & 10.97
  \end{tabular}}
  \caption*{$\rho=0.5$ and  $\eta=100$.}\label{tab:rho}
  \end{subtable}
  \caption{Results for the systemic risk problem.
  Estimations of the value function at point $(x_0=10)$ at time 0 provided by different strategies.
  We took $T=1$, $N=100$, $\sigma=0.1$, $\kappa=0.5$, $X_0=10$. }\label{t:SR}
\end{table}

\section{Conclusion}

In this work, we have investigated how to use probabilistic numerical methods for some classes of mean field control problem via Markovian embedding.
We focused on two types of Regression Monte Carlo methods (namely, Regress-Later and Control Randomization) and Quantization.
We have then presented three different examples of applications.

We found that the Regression Monte Carlo algorithms perform well in problems of control of the drift.
In such problems, they are much faster than Quantization for similar precision.
In particular, we noticed that Regress-Later is usually more reliable than Control Randomization; often the choice of a uniform distribution of the training points on an appropriate interval is sufficient to obtain high-quality estimations.
On the other hand Control Randomization is very sensitive to the choice of the distribution of the randomized control, and often a few iterations are necessary before finding a good control distribution.
We have also tried to use the performance iteration or path recomputation method, but on the examples we considered, it was very time consuming and did not help much in terms of accuracy.
Despite the success of Regression Monte Carlo methods in problems with control on the drift, the example of Portfolio Selection highlighted a possible weakness of these algorithms.
When the control acts on the variance only, we found difficult to make the numerical scheme converge to sensible results within the computational resources available.
We realized that the study of these problems and the solution via Regression Monte Carlo methods is outside the scope of this paper.
This is probably related to another limitation of this family of methods: the choice of the basis functions for the regression.
Indeed, for some problems, a good basis might be very large or might require several steps of trials and errors.

Quantization-based method, on the other hand, provided very stable and accurate results.
A first interesting and practical feature of the Q-algorithm is that regressing the value function using quantization-based methods is local.
So, first, it can be easily parallelized to provide fast results, and, second, it is easy to check at which points of the grids the estimations suffer from instability and how to change the grid to fix the problem (basically, by adding more points where the estimations need to be improved).
Another interesting feature of the quantization methods is that, one can choose the grids on which to learn the value function.
It is possible to exploit this feature in the case where one has, a priori, a rough idea of where the controlled process should be driven by the optimal strategy (see, e.g., the liquidation problem).
In this case, one should build grids with many points located where the process is supposed to go.
In the case where one has no guess of where the optimal process goes, it is always possible to use bootstrapping methods to build better grids iteratively, starting from a random guess for the grid (see, e.g., the systemic risk problem).
In both cases, one has to be particularly careful with the borders of the grids that have been built.
Indeed, the decisions computed by quantization-based methods at the borders might easily be wrong if the grids do not have a ``good shape'' at the borders.
Unfortunately, the shape of the grid that should be used depends heavily on the problem under consideration.
Except in special cases, it seems not possible to avoid the use of deterministic algorithms (such as gradient descent methods or extensive search) to find the optimal action at each point of the grid.
A smooth expression of the conditional expectations of the quantized processes is necessary for the deterministic algorithms find optimal strategy efficiently.
Once again, the use of parallel computing can alleviate the time-consuming task of searching for the optimal control at each point of the grids.

{\bf Acknowledgments.}
Most of this work was realized while the third author was a postdoctoral fellow at NYU Shanghai; the support of a research discretionary fund from the NYU-ECNU Institute of Mathematical Sciences and the support provided by the CEMRACS for his stay at CIRM are gratefully acknowledged.

\bibliographystyle{abbrv}
\bibliography{stocmkv-proceeding}

\end{document}